\documentstyle[11pt]{article} 
\begin{document} 
\newtheorem{Def}{Definition}[section] 
\newtheorem{thm}{Theorem}[section] 
\newtheorem{lem}{Lemma}[section] 
\newtheorem{rem}{Remark}[section] 
\newtheorem{prop}{Proposition}[section] 
\newtheorem{cor}{Corollary}[section] 
\title 
{Compactness of solutions to the Yamabe problem. III } 
\author{ YanYan Li\thanks{Part of the work is supported by NSF-DMS-0401118. 
Department of Mathematics, Rutgers University,110 Frelinghuysen 
Road, Piscataway, NJ 08854. email: yyli@math.rutgers.edu}
\qquad Lei Zhang \thanks{Part of the work 
is supported by NSF-DMS-0600275. Department of Mathematics, 
University of Florida, 358 Little Hall, Gainesville, FL 32611-8105.
email: leizhang@math.ufl.edu}}
\date{} 
\maketitle 
\input { amssym.def} 

\begin{abstract}
For a sequence of blow up solutions of the Yamabe equation
on non-locally confonformally flat compact Riemannian manifolds
of dimension $10$ or $11$, we establish sharp estimates on 
its  asymptotic profile near blow up points as well as sharp decay
estimates of the Weyl tensor and its covariant derivatives
at blow up points.   If the Positive Mass Theorem held
in dimensions $10$ and $11$, these estimates 
would imply the compactness of
the set of  solutions of the Yamabe equation on such
manifolds.
\end{abstract}
 
\setcounter{section}{0} 
 
\section{Introduction} 
 
Let $(M^n, g)$ be a 
compact, smooth, connected Riemannian manifold (without boundary) 
of dimension $n\ge 3$.  The Yamabe conjecture has been proved 
through the works of Yamabe \cite{Y}, Trudinger \cite{T}, 
Aubin \cite{A} 
and Schoen \cite{S}:  The conformal class of $g$ contains 
a metric of constant scalar curvature.
Different proofs of the Yamabe conjecture in the case
$n\le 5$ and in the case $(M, g)$ is locally conformally flat
are given by Bahri and Brezis \cite{BB} and Bahri \cite{B}.
 
Consider the Yamabe equation 
and its sub-critical 
approximations: 
\begin{equation} 
-\Delta_g u +c(n)R_g u 
=n(n-2)u^p,\qquad u>0,\quad \mbox{on}\quad M, 
\label{Y0} 
\end{equation} 
where $1<p\le \frac{n+2}{n-2}$, 
 $\Delta_g$ is the Laplace-Beltrami operator 
associated with $g$, $R_g$ is the scalar curvature of $g$, and 
$c(n)=\frac { (n-2) }{ 4(n-1) }$. Let 
$$ 
{\cal M}_p=\{ u\in C^2(M)\ |\ u\ 
\mbox{satisfies}\ (\ref{Y0})\}. 
$$ 
 
If $(M,g)$ is locally conformally flat and is not conformally 
diffeomorphic to the standard sphere, Schoen \cite{S2} proved  that 
for any $1<1+\epsilon\le p\le \frac {n+2}{n-2}$ and any non-negative 
integer $k$, 
\begin{equation} 
\|u\|_{ C^k(M,g) } \le C,\qquad \forall\ u\in {\cal M}_p, 
\label{compact} 
\end{equation} 
where $C$ is some constant depending only on $(M,g)$, $\epsilon$ and 
$k$. The same conclusion has been proved to hold in dimension $n\le 
7$ for 
  $(M^n,g)$ which are not locally conformally flat, see Li and Zhang 
\cite{LiZhang1} and Marques \cite{M}.
See also the introduction of \cite{LiZhang1}
where works of Li-Zhu \cite{LiZhu}, Li-Zhang \cite{LZ1} and Druet
\cite{D1}-\cite{D2} for dimensions $n=3,4,5$ are 
described.  Extensive works on the problem and closely related ones
can be found in  \cite{LiZhang1} as well.

For  $n=8, 9$ and on  $(M^n,g)$ which are not locally conformally
flat,
\cite{LiZhang1}
also contains  sharp estimates on blow up solutions
of (\ref{Y0})  and sharp decay estimates of  the Weyl
tensor and its first covariant derivatives
at  blow up points. If
 the Positive Mass Theorem held
in dimensions $8$ and $9$, these estimates would
yield (\ref{compact}) for $n=8,9$.
Soon after completing \cite{LiZhang1},  we extended
these sharp estimates
to dimensions $n=10, 11$ (see Theorem \ref{thm0});
however
 we have encountered some difficulty
in extending such  estimates to $n\ge 12$.
Very interesting  results have
subsequently been obtained  by Aubin in \cite{A1, A2}.

To study the compactness of solutions to the Yamabe equation,
it is crucial to establish  sharp
estimates of blow up solutions.
An important step is to find out
the right asymptotic profile
 of blow up solutions near a blow up point.
Our earlier work 
 \cite{LiZhang1} strongly suggested such a profile in
dimensions $n\ge 10$, which we describe below.

Let  $\{u_k\}$ be a sequence of 
solutions to the Yamabe equation on $(M^n, g)$ satisfying,
for some $\bar P_k\in M$,
$$
u_k(\bar P_k):= \max_{M} u_k\to \infty.
$$
Assume that $g=g_{ij}(z)dz^i dz^j$ is 
already in conformal normal coordinates centered
at $\bar P_k$.

For dimension $n=8,9$, sharp 
 estimates on blow up solutions
and sharp decay estimates of the  Weyl tensor and its covariant
derivatives at blow up points
were established in  \cite{LiZhang1}
through an iterative procedure. Due to  this
 procedure,
 we expect to obtain enough
estimates on the decay rates of the Weyl tensor
and its covariant derivatives of appropriate
order before making the next step in the
iterative  process, and  therefore we can
use \begin{equation}
-\Delta u_k +c(n)R_g u_k = n(n-2)
u_k^{ \frac {n+2}{n-2}},
\label{app1}
\end{equation}
instead of the Yamabe equation which would 
replace  $\Delta $ in  (\ref{app1})
by $\Delta_g$, 
to determine the  asymptotic profile of $\{u_k\}$
near blow up points.
Note that  $\Delta$ 
 is the flat Laplacian in the $z-$coordinates.

The Taylor expansion of $R(z)$
in conformal normal coordinates  is, for  $\bar l\ge 2$,
\begin{equation}
R(z)=
\sum_{l=2}^{\bar l} \sum_{|\alpha|=l} \frac {\partial_\alpha R}{\alpha!}
z^\alpha+O(|z|^{\bar l+1}).
\label{R1}
\end{equation}
For convenience, we write, with
$M_k:= u_k(0)$, 
$$
v_k(y):=M_k^{-1}u_k(M_k^{-\frac{2}{n-2}}y).
$$
Then (\ref{app1}) becomes
\begin{equation}
\Delta v_k(y)-\bar c(y) 
v_k(y)+n(n-2)v_k(y)^{\frac{n+2}{n-2}}=0
\label{app3}
\end{equation}
where
$$
\bar c(y):= c(n)R_g(M_k^{-\frac 2{n-2}}y)M_k^{-\frac 4{n-2}}
=c(n) \sum_{l=2}^{\bar l} 
M_k^{ -\frac{4+2l}{n-2} }  \sum_{|\alpha|=l} \frac {\partial_\alpha R}{\alpha!}
y^\alpha
$$
and
$$
v_k \ \mbox{converges }\ 
\mbox{to}\ U(y):=(\frac 1{1+|y|^2})^{ \frac {n-2}2 }
\ \mbox{in}\ C^3_{loc}(\Bbb R^n).
$$
For dimension $n=10,11$, we only need to consider, in the
formal expansion of $v_k$, 
$$
\tilde v_k= v^{(1)}+ M_k^{ -\frac 8{n-2} }
v^{(2)}+ M_k ^{ -\frac {10}{n-2} }
v^{(3)}.
$$
The equations satisfied by $v^{(1)}, v^{(2)}$ and $v^{(3)}$ are,
determined by (\ref{app3}),
$$
\Delta  v^{(1)}+ n(n-2) [v^{(1)}]^{ \frac {n+2}{n-2} }=0,
$$
$$
\Delta  v^{(2)} - c(n)[ \sum_{|\alpha|=2} \frac {\partial_\alpha R}{\alpha!}
y^\alpha]
v^{(1)} +n(n+2) [v^{(1)}]^{ \frac 4{n-2} }  v^{(2)}=0,
$$
and
$$
\Delta  v^{(3)} - c(n)[ \sum_{|\alpha|=3} \frac {\partial_\alpha R}{\alpha!}
y^\alpha]
v^{(1)} +n(n+2) [v^{(1)}]^{ \frac 4{n-2} }  v^{(3)}=0.
$$

Let
\begin{equation}
\bar R^{(l)}:=\frac 1{ |\Bbb S^{n-1}| }
\int_{\theta\in \Bbb S^{n-1} }
\sum_{|\alpha|=l} \frac {\partial_\alpha R}{\alpha!}
\theta^\alpha,
\label{R2}
\end{equation}
and
\begin{equation}
\tilde  R^{(l)}(\theta):= -\bar R^{(l)}+
\sum_{|\alpha|=l} \frac {\partial_\alpha R}{\alpha!}
\theta^\alpha, \qquad\theta\in \Bbb S^{n-1}.
\label{R3}
\end{equation}
We have, in polar coordinates, 
\begin{equation}
 v^{(1)}:= U,
\quad v^{(2)}(r, \theta)
= -c(n)\tilde R^{(2)}(\theta)f_2(r),
\quad v^{(3)} (r, \theta)
= -c(n)\tilde R^{(3)}(\theta)f_3(r),
\label{V9}
\end{equation}
where $f_2$ and $f_3$ 
are respectively the solutions of (\ref{nov19e4}) and (\ref{aabb}).

Thus 
the  asymptotic profile of $u_k$ near the blow up point
$\bar P_k$  should be 
\begin{eqnarray}
\tilde u_k(x)&=& M_k \tilde v_k(M_k ^{\frac 2{n-2} }x)
\nonumber\\
&=& u_k(0)\big[ v^{(1)} (u_k(0)^{ \frac 2{n-2} }x)+
u_k(0)^{ -\frac 8{n-2} } v^{(2)}(u_k(0)^{  \frac 2{n-2} }x)+
 u_k(0)^{ -\frac {10}{n-2} } v^{(3)}(u_k(0)^{  \frac 2{n-2} }x)\big]\nonumber\\
&=& u_k(0) U(u_k(0) ^{\frac 2{n-2} }x)
-
 u_k(0) ^{ \frac{n-10}{n-2} }
\big[c(n) \tilde R^{(2)}(\theta) f_2(u_k(0) ^{\frac 2{n-2} }|x|)\big]
\nonumber\\
&& - u_k(0) ^{ -\frac{12-n}{n-2} }
\big[c(n) \tilde R^{(3)}(\theta) f_3(u_k(0) ^{\frac 2{n-2} }|x|)\big].
\nonumber
\end{eqnarray}

In the  following, we give  
 the previously mentioned 
 sharp estimates
in dimensions  $n=10,11$.  The asymptotic  profile of blow up
solutions 
is exactly the one described above.   

For $Q\in M$ and $\mu>0$, let
$$\xi_{Q,\mu}(P):=
(\frac{\mu}{1+\mu^2\mbox{dist}_g(P,Q)^2})^{\frac{n-2}{2}},
\quad P\in M 
$$
and, in polar coordinates, 
\begin{eqnarray*}
\tilde \xi_{Q,\mu}(P)&=&
\xi_{Q,\mu}(P)-c(n)\tilde R^{(2)}(\theta)
f_2(\mu\cdot \mbox{dist}_g(P,Q))\mu^{\frac{n-10}{2}}
\\
&&-c(n)\tilde R^{(3)}(\theta)
f_3(\mu\cdot \mbox{dist}_g(P,Q))\mu^{\frac{n-12}2}.
\end{eqnarray*}

We use $W_g$ to denote  the Weyl tensor of the metric $g$.

\begin{thm}  Let $(M^n,g)$ be a compact, smooth, connected 
Riemannian manifold of dimension $n=10,11$, and let
 $u$ be a smooth  solution of (\ref{Y0})
with  $1<1+\epsilon\le p\le
\frac{n+2}{n-2}$. 
  Then for  some 
positive constant 
$C$ and some positive integer $m$ which 
 depend only 
 on 
  $(M,g)$, 
  there exist  some local maximum points of $u$,
  denoted as 
${\cal S}:=\{P_1, \cdots, P_m\}$, 
such that 
$$ \mbox{dist}_g(P_i,P_j)\ge \frac 1{C},\qquad
\frac 1Cu(P_i)\le u(P_j)\le Cu(P_i),\qquad  \forall 
i\neq j, $$ 
$$ 
|W_g(P_i)|_g 
\le C u(P_i)^{  -\frac{n-6}{n-2}  }, 
 \qquad |\nabla W_g(P_i)|_g\le  C u(P_i)^{  -\frac{n-8}{n-2}  }, 
$$ 
 \begin{eqnarray*} 
 |\nabla_g^2  W_g(P_i)|_g&\le \left\{ 
 \begin{array}{ll} 
\frac C{ \sqrt{  \log u(P_i) } },,& \mbox{if}\ n=10,\\ 
C 
 u(P_i)^{  -\frac{n-10}{n-2}  },& \mbox{if}\ n=11, 
 \end{array} 
 \right. 
\qquad\quad \forall\ i,
 \end{eqnarray*} 
$$
\frac 1C\sum_{l=1}^{m}\xi_{P_l,u(P_l)^{ \frac 2{n-2} }}\le u\le
C\sum_{l=1}^{m}\xi_{P_l,u(P_l)^{ \frac 2{n-2} }},
\quad \ \mbox{on}\ M.
$$
Moreover, for each $l$ and  modulo a conformal factor
 which makes $g$ in  conformal normal coordinates
at $P_l$,
\begin{eqnarray*}
&&|\nabla_g^{\alpha}(u-
\tilde \xi_{P_l,u(P_l)^{\frac 2{n-2}}})(P)|\\&\le
 &Cu(\bar P)^{\frac{n-14+2|\alpha|}{n-2}}(1+u(\bar P)^{\frac 2{n-2}}
\mbox{dist}_g(P,P_l)^{8-n-|\alpha |}),
\quad \forall\ dist_g(P, P_l)<\frac 1C,\
 |\alpha |=0,1,2.
\end{eqnarray*}
\label{thm0} 
\end{thm}

A consequence of Theorem \ref{thm0} is 
\begin{cor} 
 Let $(M^n,g)$, $n=10,11$, be a 
compact, smooth, connected Riemannian manifold  which is not locally 
conformally flat, and let
  $1<1+\epsilon\le p\le 
\frac{n+2}{n-2}$.  Then 
$$ 
\|u\|_{ H^1(M,g) } \le C,\qquad \forall\ u\in {\cal M}_p, 
$$ 
where  $C$ is some constant depending only on 
$(M^n,g)$ and  $\epsilon$. 
\label{cor1} 
\end{cor}

\begin{rem} If the
Positive Mass Theorem held in 
 dimensions $n=10$ and $11$, 
 Theorem \ref{thm0} would yield 
 (\ref{compact}) in these dimensions.
\label{rem1} 
\end{rem}

In the following we give a result which
is more local in nature.
Let $B_1\subset \Bbb R^n$   be the 
unit ball centered at the origin, and let 
$(a_{ij}(x))$ be a smooth,  $n\times n$ symmetric 
positive definite matrix function, defined 
on $B_1$,  satisfying 
\begin{equation} 
\frac 12 |\xi|^2\le a_{ij}(x)\xi^i\xi^j\le 2  |\xi|^2, 
\qquad \forall\ x\in B_1, \ \xi\in \Bbb R^n, 
\label{cond1} 
\end{equation} 
and, for some $\bar a>0$, 
\begin{equation} 
\|a_{ij}\|_{ C^{5}(B_1) }\le \bar a. 
\label{cond2} 
\end{equation} 

Consider 
\begin{equation}
-L_g u=n(n-2)u^p,\qquad u>0,\quad \mbox{on}\quad B_1,
\label{eq3}
\end{equation}
where
\begin{equation}
g:=a_{ij}(x)dx^idx^j.
\label{2-1}
\end{equation}
If $\{x^1, \cdots, x^n\}$  are conformal normal coordinates
for $g$, let
\begin{equation}
v:= v^{(1)}+ u(0)^{ -\frac 8{n-2} }v^{(2)}+
u(0)^{ -\frac{10}{n-2} }v^{(3)},
\label{V7}
\end{equation}
where $ v^{(1)}, v^{(2)}$ and $v^{(3)}$ are defined in (\ref{V9}).
 
\begin{thm} 
 Let $(B_1,g)$ be as  above 
 and let 
 $u$ be  a solution of (\ref{eq3}), with $1<1+ \epsilon 
\le p\le \frac {n+2}{n-2}$ and $n=10,11$.
Assume,
  for some constant $\bar b\ge 1$, 
\begin{equation} 
\nabla u(0)=0, \qquad 1\le \sup_{ B_1}u\le \bar b u(0). 
\label{normal} 
\end{equation} 
 Then  there exist some 
positive constants $\delta$ and $C$, depending only on 
$\bar b$, $\epsilon$ and $\bar a$,  such 
that 
\begin{equation} 
u(0)u(x)|x|^{n-2}\le 
C, \qquad \forall\ 0<|x|\le\delta, 
\label{eq4} 
\end{equation}
 
\begin{eqnarray} 
&&|W_g(0)|_g ^2 u(0)^{-\frac 8{n-2}} 
+|\nabla_gW_g(0)|_g^2  u(0)^{-\frac {12}{n-2}} 
+|\nabla_g ^2W_g(0)|_g ^2  u(0)^{-\frac {16}{n-2}} 
\log u(0)\nonumber\\ 
&\le& C u(0)^{-2},
\qquad\qquad\qquad\qquad\qquad\qquad\qquad\qquad\quad n=10, 
\label{W1} 
\end{eqnarray} 
\begin{eqnarray} 
&&|W_g(0)|_g ^2 u(0)^{-\frac 8{n-2}} 
+|\nabla_gW_g(0)|_g ^2 u(0)^{-\frac {12}{n-2}} 
+|\nabla_g ^2W_g(0)|_g ^2 u(0)^{-\frac {16}{n-2}}\nonumber\\ 
&\le& C u(0)^{-2}, \qquad\qquad\qquad\qquad\qquad\qquad\qquad\quad  n=11,
\label{W2} 
\end{eqnarray}
\begin{equation}
\frac 1{C} u(0)U(u(0)^{\frac 2{n-2}}x)
\le u(x)\le C u(0)U(u(0)^{\frac 2{n-2}}x),
 \qquad \forall\ 0<|x|\le\delta,
\label{W7}
\end{equation}
and, if
$\{x^1, \cdots, x^n\}$  are conformal normal coordinates
for $g$,  we have, with $v$ given in (\ref{V7}),  
\begin{eqnarray}
&&|\nabla^\alpha(u- u(0) v(u(0)^{\frac 2{n-2}}\cdot))|
\nonumber\\
&\le& C
u(0)^{\frac{n-14+2|\alpha|}{n-2}}(1+u(0)^{\frac 2{n-2}}
|x|)^{8-n-|\alpha |}),
 \qquad \forall\ 0<|x|\le\delta, \
|\alpha|=0, 1, 2.
\label{W8}
\end{eqnarray}
\label{thm1} 
\end{thm}

It is not difficult to see that
 Theorem \ref{thm0} follows from Theorem \ref{thm1}.
Our proof of Theorem \ref{thm1} follows  closely the 
arguments in \cite{LiZhang1}.  
In particular the sharp 
estimates on blow up solutions
and on the decay rates of the  Weyl tensor and its  derivatives
 are
 obtained iteratively with improved estimates after each 
iteration.
The main difference between the arguments in this paper and those of 
\cite{LiZhang1} is that some Riemannian tensor inequalities
in  
 conformal normal coordinates which we used for dimension $n=8,9$ 
are not sufficient for higher dimensions.
Our proof of  Theorem \ref{thm1}
requires an estimate from below of some integral
quantity associated with $v^{(2)}$.

We will only prove Theorem \ref{thm1} for $p=\frac{n+2}{n-2}$ since
modifications of the arguments yield the result for
$1+\epsilon\le p\le \frac{n+2}{n-2}$, see section 5 of
\cite{LiZhang1}. We will always assume that $n=10,11$
unless otherwise stated.
 In Section 2 we 
prove Theorem \ref{thm1}. 
In the appendices we establish some facts 
which we use for the proof.

\section{The  proof of Theorem \ref{thm1}} 
In this section we prove  Theorem \ref{thm1}.
In the first four subsections we establish 
(\ref{eq4}) using  the method of moving spheres.
In the last section we 
derive (\ref{W1})-(\ref{W8})  using the
Pohozaev type identity.

\subsection{The set up for proving  (\ref{eq4})} 
  Suppose the contrary of 
 (\ref{eq4}), 
then for some $\bar a>0$, $\bar b\ge 1$, there exist a sequence of 
Riemannian metrics $\{\tilde g_k\}$ of the form  (\ref{2-1}) that 
satisfy (\ref{cond1}) and (\ref{cond2}), and some solutions $u_k$ of 
(\ref{eq3}), with $p=\frac {n+2}{n-2}$ and with $g$ replaced by 
$\tilde g_k$, satisfying (\ref{normal}), such that 
\begin{equation} 
\max_{ |x|<\frac 1k} \bigg( 
u_k(0) u_k(x) |x|^{n-2} \bigg)\ge k, 
\label{5-0} 
\end{equation} 
 
We will simply use $g$ to denote $\tilde g_k$, and we assume that 
$g_{ij}(z)dz^idz^j$ is  already in conformal normal coordinates 
centered at the origin --- as in the proof of theorem 2.1 
 in  \cite{LiZhang1}. 
As in  \cite{LiZhang1}, 
$$
M_k:=u_k(0)\to \infty.
$$ 
Write 
$$
(g_k)_{ij}(y)=g_{ij}(M_k^{ -\frac 2{n-2} }y)dy^idy^j,
$$ 
$$ 
v_k(y):=M_k^{-1}u_k(M_k^{-\frac{2}{n-2}}y), 
$$ 
$$ 
c(x)=c(n)R_g(x), \quad \mbox{and} \quad 
 \bar c(y)=c(n)R_g(M_k^{-\frac 2{n-2}}y)M_k^{-\frac 4{n-2}}. 
$$ 
Then 
\begin{equation} 
\left\{ 
\begin{array}{ll} 
\Delta_{g_k}v_k(y)-\bar cv_k(y)+n(n-2)v_k(y)^{\frac{n+2}{n-2}}=0, 
\quad |y|\le \frac 12M_k^{\frac{2}{n-2}}, 
  \\ 
1=v_k(0)\ge (\bar b^{-1}+\circ(1))v_k(y), 
\quad |y|\le \frac 12M_k^{\frac{2}{n-2}},\quad 
\nabla v_k(0)=0. 
\end{array} 
\right. 
\label{ab1} 
\end{equation} 
By the Liouville type theorem of Caffarelli, Gidas and Spruck \cite{CGS},
together with some standard elliptic estimates, 
$$
v_k \ \mbox{converges}\ \mbox{to}\ U(y):=(\frac 1{1+|y|^2})^{ \frac {n-2}2 }
\ \mbox{in}\ C^3_{loc}(\Bbb R^n).
$$
In local coordinates, 
\begin{eqnarray} 
g_{pq}(x)&=&\delta_{pq}+\frac 13R_{pijq}x^ix^j+\frac 16R_{pijq,k}x^ix^jx^k 
\nonumber \\ 
&&+(\frac 1{20}R_{pijq,kl}+\frac 2{45}R_{pijm}R_{qklm})x^ix^jx^kx^l+O(r^5). 
\nonumber 
\end{eqnarray} 
In conformal normal coordinates, write 
$$ 
\Delta_g=\frac 1{\sqrt{g}}\partial_i(\sqrt{g}g^{ij}\partial_j) 
=\Delta +b_i\partial_i +d_{ij}\partial_{ij}, 
$$ 
where $(g^{ij})$ denotes the inverse matrix of $(g_{ij})$, 
$\partial_i=\frac {\partial }{\partial z^i}$, 
 $\partial_{ij}=\frac {\partial ^2}{ \partial z^i \partial z^j}$, 
$\Delta=\sum_{i=1}^n \frac {\partial ^2}{ \partial z^i \partial z^i}$, 
\begin{eqnarray} 
&&b_i(x)= \partial_j g^{ij}(x)\nonumber\\ 
&=&-\frac 16 R_{ia,b}x^ax^b-\frac 16R_{iabp,p}x^ax^b-(\frac 1{20} 
R_{ia,bc}-\frac 1{15}R_{ipad}R_{pbcd} 
\nonumber\\ 
&&-\frac 1{15}R_{iapd}R_{pbcd}+\frac 
 1{10}R_{iabp,pc})x^ax^bx^c+O(r^4), \nonumber 
\end{eqnarray} 
and 
\begin{eqnarray} 
d_{ij}(x)=g^{ij}-\delta_{ij}=-\frac 13R_{ipqj}x^px^q-\frac 16R_{ipqj,k}x^px^qx^k-(\frac 1{20}R_{ipqj,kl}\nonumber \\ 
-\frac 1{15}R_{ipqm}R_{jklm})x^px^qx^kx^l+O(r^5).\nonumber 
\end{eqnarray} 
Thus 
$$\Delta_{g_k}=\Delta+\bar b_i\partial_i+\bar d_{ij}\partial_{ij}, 
$$ 
where 
$$ 
 \bar b_i(y)=M_k^{-\frac{2}{n-2}}b_i(M_k^{-\frac{2}{n-2}}y), 
\quad \bar d_{ij}(y)=d_{ij}(M_k^{-\frac{2}{n-2}}y). 
$$ 
 
For  $\lambda>0$ and for any function $v$, let, as in 
\cite{LiZhang1}, 
$$v^{\lambda}(y):=(\frac{\lambda}{|y|})^{n-2} 
v(y^\lambda),  \qquad y^\lambda:=\frac{\lambda^2y}{|y|^2},$$ 
denote the Kelvin transformation of $v$, and 
$$ 
\Sigma_{\lambda}:=B(0, \frac 1{\sqrt{k}} 
M_k^{\frac{2}{n-2}})\setminus \overline{B(0, \lambda)}=\{y\ |\ 
\lambda<|y|< 
 \frac 1{\sqrt k}  M_k^{\frac{2}{n-2}}\}, 
$$ 
$$ 
w_{\lambda}(y):= v_k(y)-v_k^{\lambda}(y),\qquad y\in \Sigma_\lambda. 
$$ 
 
As in (33)-(35) in \cite{LiZhang1}, 
 the equation 
for $w_{\lambda}$ is 
$$ 
\Delta w_{\lambda}+\bar b_i\partial_i w_{\lambda} 
+\bar d_{ij}\partial_{ij}w_{\lambda} 
-\bar cw_{\lambda}+n(n+2)\xi^{\frac 4{n-2}}w_{\lambda}= 
E_{\lambda}, \qquad\mbox{in}\ \Sigma_\lambda, 
$$ 
where $ \xi^{ \frac 4{n-2}}=\int_0^1(tv_k+(1-t)v_k^{\lambda})^{\frac 
4{n-2}}dt$ and 
\begin{eqnarray} 
E_{\lambda}&=& \left(\bar c(y)v_k^{\lambda}(y)- (\frac{\lambda}{|y|})^{n+2} 
\bar c(y^{\lambda})v_k(y^{\lambda})\right) 
-(\bar b_i\partial_i v_k^{\lambda}+\bar d_{ij}\partial_{ij}v^{\lambda}_k) 
\nonumber\\ 
&& +(\frac{\lambda}{|y|})^{n+2}\left(\bar b_i(y^{\lambda}) 
\partial_i v_k(y^{\lambda})+ 
\bar d_{ij}(y^{\lambda})\partial_{ij}v_k(y^{\lambda})\right). 
\label{eQ} 
\end{eqnarray} 
 
 As in \cite{LiZhang1} 
 we apply the method of moving spheres to 
$w_{\lambda}+h_{\lambda}$, an appropriate 
perturbation of $w_{\lambda}$, 
 for $\lambda$ in a fixed 
 neighborhood of $1$, to reach
 a contradiction. 
We recall that the perturbation 
 $h_{\lambda}$ is required to satisfy the 
following properties: 
 
\begin{enumerate} 
\item $h_{\lambda}=0$ on $\partial B_{\lambda}$,
 $h_{\lambda}=\circ(1)|y|^{2-n}$ in
$\Sigma_\lambda$.
\item
For $O_{\lambda}:=
\{y\in \Sigma_\lambda\ |\ v_k(y)<2v_k^{\lambda}(y)\}$,
\begin{equation} 
\label{20-1new} (\Delta +\bar b_i\partial_i+ \bar 
d_{ij}\partial_{ij}-\bar c+n(n+2)\xi^{\frac{4}{n-2}})h_{\lambda} 
+E_{\lambda}\le 0\ \mbox{ in }\ O_{\lambda}.
\end{equation} 
\end{enumerate} 
Note that by the first property, $w_{\lambda}+h_{\lambda}>0$ in 
$\Sigma_{\lambda}\setminus \bar O_{\lambda}$, so we do not need 
(\ref{20-1new}) to hold outside $O_{\lambda}$.

\subsection{Estimate of  $v_k$ in $|y|\le 
M_k^{ \frac{16-\epsilon}{  (n-2)^2 } }$} 
 
In \cite{LiZhang1} we established estimates on $v_k-U$ with
an error term of the order
 $M_k^{-\frac 8{n-2}}$.
Now, for dimension 
$n=10,11$,  we need to work with terms in the
formal expansion of $v_k$ as 
described  in the introduction which are of order
$M_k^{-\frac 8{n-2}}$ and $M_k^{-\frac {10}{n-2}}$.
The main result of this subsection  
is the following estimate of $v_k$ with an error term of
the order $M_k^{-\frac {12}{n-2}}$ in the region $|y|\le M_k^{
\frac{16-\epsilon}{ (n-2)^2 } }$ for any $\epsilon>0$.
\begin{prop} For $n\ge 10$, and for any $\epsilon>0$,
there exists some positive constant $C(\epsilon)$ such that
\begin{eqnarray}
&&\bigg|\nabla^l\left(
v_k-( v^{(1)}+ M_k^{ -\frac 8{n-2} }
v^{(2)}+ M_k ^{ -\frac {10}{n-2} }
v^{(3)})\right)\bigg|
\nonumber\\
&\le& C(\epsilon)M_k^{-\frac {12}{n-2}}(1+r)^{8-n+\bar a-l},\quad
\qquad \qquad 0<r\le M_k^{\frac{16-\epsilon}{(n-2)^2}},\
l=0,1,2,
\label{nov19e6}
\end{eqnarray}
where $\bar a=\frac 34(n-10+\sqrt{\epsilon})$, $ v^{(1)}$,
$v^{(2)}$ and $v^{(3)}$ are defined in (\ref{V9}).
\label{propmar201}
\end{prop}

We first recall some notations 
in  \cite{LiZhang1}. 
For $\bar l\ge 2$, write the Taylor expansion of $R(x)$ at $0$
as (\ref{R1}).
Let $\bar R^{(l)}$ and $\tilde  R^{(l)}(\theta)$ be defined as
in (\ref{R2}) and (\ref{R3}).
We know (see (44) in \cite{LiZhang1}), with $W$ denoting
the Weyl tensor,  that 
\begin{equation} 
\bar R^{(2)}=\frac 1{2n}\Delta R= 
-\frac 1{12n}|W|^2, 
\qquad \mbox{and}\ 
\bar R^{(3)}=0. 
\label{SS5} 
\end{equation} 
 
We write 
\begin{equation} 
\tilde R^{(l)}(\theta)=\sum_ 
{p\ge 1}\tilde R_{lp}e_{p}(\theta)\quad 2\le l\le \bar l. 
\label{15-4new} 
\end{equation} 
where $e_{p}'s$, 
depending only on $n$, 
 are non-constant eigenfunctions of $-\Delta_{\Bbb S^{n-1}}$. 
The following lemma, whose proof can be found in Appendix A,
is used  in our arguments.
\begin{lem} 
\label{lemjune27} 
\begin{equation} 
\tilde R^{(3)}=\sum_{p=1}^{l_3}\tilde R^{(3)}_{p}e_{3p}(\theta)+O(|W|), 
\label{7-1new} 
\end{equation} 
where $\left\{e_{3p}(\theta)\right\}_{1\le p\le l_3}$ is a set of 
eigenfunctions of $-\Delta_{\Bbb S^{n-1} }$ associated with the 
eigenvalue $3(n+1)$. 
\end{lem} 
 
Let $f_2$ and $f_3$ be defined as  in
Appendix C,
set 
$$ 
F^{(2)}:=-c(n)\tilde R^{(2)}(\theta)f_2(r)M_k^{-\frac 8{n-2}} 
=v^{(2)} M_k^{-\frac 8{n-2}}
$$ 
and
\begin{eqnarray*}
F^{(3)}&:=&F^{(2)}-c(n)\sum_{p=1}^
{l_3}\tilde R^{(3)}_{p}e_{3p}(\theta)f_3(r)M_k^{-\frac{10}{n-2}}
\\
&=& v^{(2)} M_k^{-\frac 8{n-2}}
+  v^{(3)} M_k^{-\frac{10}{n-2}}
+O(|W|)f_3(r)M_k^{-\frac{10}{n-2}}.
\end{eqnarray*}
By
  (\ref{nov19e4}) and (\ref{aabb}), 
$$ 
(\Delta +n(n+2)U^{\frac 4{n-2}})F^{(2)}=
c(n)\tilde R^{(2)}(\theta)r^2UM_k^{-\frac 8{n-2}}, 
$$ 
and
\begin{eqnarray} 
&&(\Delta +n(n+2)U^{\frac 4{n-2}})F^{(3)} \nonumber \\ 
&=&c(n)\sum_{l=2}^3\tilde R^{(l)}(\theta)r^lUM_k^{-\frac{4+2l}{n-2}} 
+O(|W|)M_k^{-\frac {10}{n-2}}(1+r)^{5-n}. 
\label{dec7e8} 
\end{eqnarray}

\noindent{\bf Proof of Propostion \ref{propmar201}.}\ We claim that 
\begin{equation} 
(\Delta_{g_k}-\bar c)(U+F^{(3)})+n(n-2)(U+F^{(3)})^{\frac{n+2}{n-2}} 
=O(M_k^{-\frac{12}{n-2}})(1+r)^{6-n}. 
\label{dec8e1} 
\end{equation} 
To see this, we first recall some known facts. We know from (21), 
(44) and (123) in \cite{LiZhang1} that 
\begin{equation} 
\label{jul3e1} 
 v_k(y)\le CU(y),\quad |y|\le 
M_k^{\frac{16-\epsilon}{(n-2)^2}},\quad n\ge 10. 
\end{equation} 
\begin{equation} 
|\nabla^lR_{abcd}|=O(M_k^{-\frac{2(2-l)}{n-2}+\epsilon}), 
\quad l=0,1,\quad n\ge 10. 
\label{8-2} 
\end{equation} 
$$ 
\bar R^{(2s+2)}=O(M_k^{-\frac{4(2-s)}{n-2}+\epsilon}),\quad s=0,1,\quad 
n\ge 10. 
$$ 
These  lead to 
\begin{eqnarray} 
\bar b_i(y)&=&O(M_k^{-\frac{8-\epsilon}{n-2}})r^2+O(M_k^{-\frac{8}{n-2}})r^3, \nonumber \\ 
\bar d_{ij}(y)&=&O(M_k^{-\frac{8-\epsilon}{n-2}})(1+r)^3+O(M_k^{-\frac{8}{n-2}})r^4. 
\label{dec9e1} 
\end{eqnarray} 
 
To derive (\ref{dec8e1}) we use (\ref{dec7e8}) and (\ref{8-2}) to 
obtain 
\begin{eqnarray*} 
&&\Delta (U+F^{(3)})+n(n-2)(U+F^{(3)})^{\frac{n+2}{n-2}}\\ 
&=&\Delta (U+F^{(3)})+n(n-2)U^{\frac{n+2}{n-2}} 
(1+\frac{F^{(3)}}{U})^{\frac{n+2}{n-2}}\\ 
&=&\Delta U+n(n-2)U^{\frac{n+2}{n-2}}+ 
\Delta F^{(3)}+n(n+2)U^{\frac{4}{n-2}}F^{(3)}+ 
O(\frac{F^{(3)}}{U})^2U^{\frac{n+2}{n-2}}\\ 
&=&c(n)\sum_{l=2}^3\tilde R^{(l)}r^lUM_k^{-\frac{4+2l}{n-2}} 
+O(|W|)M_k^{-\frac{10}{n-2}}(1+r)^{5-n} 
+O(M_k^{-\frac{16}{n-2}})(1+r)^{6-n}\\ 
&=&c(n)\sum_{l=2}^3\tilde R^{(l)}r^lUM_k^{-\frac{4+2l}{n-2}} 
+O(M_k^{-\frac{14-\epsilon}{n-2}})(1+r)^{5-n}. 
\end{eqnarray*} 
Note that we used the estimates of $f_2$ and $f_3$ in 
Appendix C.
To estimate $(\Delta_{g_k}-\Delta)(U+F^{(3)})$, we observe that for 
any smooth functions $a(\theta)$ and $b(r)$,  by the definition of 
the conformal normal coordinates, $(\Delta_{g_k}-\Delta)b(r)=0$, 
consequently 
\begin{equation} 
(\Delta_{g_k}-\Delta)\bigg (a(\theta)b(r)\bigg )=\bigg ((\Delta_{g_k}-\Delta)a(\theta) 
\bigg )b(r). 
\label{925e1} 
\end{equation}

It follows, 
using the estimates of $\bar b_i$ and $\bar d_{ij}$ in (\ref{dec9e1}), 
that 
\begin{eqnarray*} 
&&(\bar b_i\partial_i+\bar d_{ij}\partial_{ij})(U+F^{(3)}) 
=(\bar b_i\partial_i+\bar d_{ij}\partial_{ij})F^{(3)}\\ 
&=&O(M_k^{-\frac{16-\epsilon}{n-2}})(1+r)^{7-n} 
+O(M_k^{-\frac{16}{n-2}})(1+r)^{8-n}. 
\end{eqnarray*} 
 
Also, by $\bar R^{(2)}=O(M_k^{-\frac{8-\epsilon}{n-2}})$ and the estimates of 
$f_2$ and $f_3$, 
$$\bar c(U+F^{(3)})=c(n)\sum_{l=2}^3\tilde R^{(l)}r^lUM_k^{-\frac{4+2l}{n-2}} 
+O(M_k^{-\frac{12}{n-2}})(1+r)^{6-n}.$$ 
 
Then (\ref{dec8e1}) is the consequence of the above. 

By (\ref{dec8e1}) and the equation for $v_k$, we have 
\begin{eqnarray} 
&&(\Delta_{g_k}-\bar c)(v_k-U-F^{(3)}) 
+n(n-2)(v_k^{\frac{n+2}{n-2}}-(U+F^{(3)})^{\frac{n+2}{n-2}}) 
\nonumber\\ 
&=&O(M_k^{-\frac{12}{n-2}})(1+r)^{6-n},\quad |y|\le 
M_k^{\frac{16-\epsilon}{(n-2)^2}}. 
\label{9-0} 
\end{eqnarray} 
Since  $f_2'(0)=f'_3(0)=0$,  $\nabla 
(U+F^{(3)})(0)=0$. By these facts and (\ref{jul3e1}) we can prove 
\begin{equation} 
\Lambda_k:= 
\max_{ |y|\le M_k^{\frac{16-\epsilon}{(n-2)^2}}} 
 |(v_k-U-F^{(3)})(y)|\le CM_k^{-\frac{12}{n-2}}. 
\label{9-1} 
\end{equation} 
Indeed, let 
$$ 
w_k:= \Lambda_k^{-1} (v_k-U-F^{(3)}). 
$$ 
Then we see from (\ref{9-0}) and (\ref{8-2}) that, 
for some $\bar\epsilon>0$ independent of $k$, 
\begin{eqnarray*} 
&&\left( \Delta+ \frac{  \circ(1) \partial_{ij}} 
{  (1+|y|)^{\bar\epsilon}  }+ 
 \frac{  \circ(1) \partial_i} 
{  (1+|y|)^{1+\bar\epsilon}  }+ \frac{  \circ(1)} 
{  (1+|y|)^{2+\bar\epsilon}  }\right)w_k(y)\\ 
&=& O(1) \Lambda_k^{-1} M_k^{ -\frac{12}{n-2}  } 
 (1+|y|)^{6-n}  +O(1)  (1+|y|)^{-4}w_k\\ 
&=& O(1) \Lambda_k^{-1} M_k^{ -\frac{12}{n-2}  } 
 (1+|y|)^{-2-\bar\epsilon} + 
O(1)  (1+|y|)^{-2-\bar\epsilon},\qquad 
|y|\le M_k^{  \frac{16-\epsilon}{(n-2)^2}  }. 
\end{eqnarray*} 
If (\ref{9-1}) did  not hold, then $ \Lambda_k^{-1} M_k^{ 
-\frac{12}{n-2}  } =\circ(1)$ along a subsequence, and the argument 
below (101) in \cite{LiZhang1} (with $\delta R_k$ replaced by $M_k^{ 
\frac{16-\epsilon}{(n-2)^2}  }$) yields a contradiction.  See also 
lemma 3.3 in \cite{ChenLin2} for a similar argument. (\ref{nov19e6}) 
is proved for $l=0$ and $|y|<R$ for $R$ being a fixed large 
constant. Next we use (\ref{9-1}) to compare $(v_k-U-F^{(3)})(y)$ 
with $QM_k^{-\frac{12}{n-2}}r^{8-n+\bar a}$ for some large $Q$ over 
$R<|y|<M_k^{\frac{16-\epsilon}{(n-2)^2}}$. By the maximum principle, 
$$ 
|v_k-(U+F^{(3)})|\le QM_k^{-\frac{12}{n-2}}r^{8-n+\bar a}. 
$$ 
The estimates for the first and the second derivatives of 
$v_k-(U+F^{(3)})$ follow from this and the equation for 
$v_k-(U+F^{(3)})$ by elliptic estimates. Proposition 
\ref{propmar201} is established. $\Box$

\subsection{Estimate of $E_{\lambda}$} 
 
In this and the next 
 subsections, we assume $\lambda\in (\frac 12, 2)$ and we 
 assume 
$\lambda\le |y|\le \frac 12 
M_k^{ \frac2 {n-2} }$ 
 unless otherwise stated. 
We use $E_1, \cdots, E_4$ to denote the following terms: 
\begin{eqnarray*} 
&&E_1=c(n)U^{\lambda}\sum_{s=0}^{2}\bar 
R^{(2s+2)}M_k^{-\frac{8+4s}{n-2}} 
r^{2s+2}(1-(\frac{\lambda}{r})^{4s+8})  \\ 
&&-\frac{c(n)^2}{2n(n+2)}\bigg (\sum_{i<j} 
2(\partial_{ij}R)^2+\sum_{i}(\partial_{ii}R)^2\bigg 
)(1-(\frac{\lambda}{r})^{8})r^2f_2^{\lambda} M_k^{-\frac{16}{n-2}} 
\end{eqnarray*} 
 
$$ 
E_2=\left\{\begin{array}{ll} 
c(n)U^{\lambda}\sum_{l=2}^{6}\tilde R^{(l)} 
M_k^{-\frac{4+2l}{n-2}}r^l(1-(\frac{\lambda}{r})^{2l+4}),\quad n=10,\\ 
\\ 
c(n)U^{\lambda}\sum_{l=2}^{7}\tilde R^{(l)} 
M_k^{-\frac{4+2l}{n-2}}r^l(1-(\frac{\lambda}{r})^{2l+4}),\quad n\ge11. 
\end{array} 
\right. 
$$

\begin{equation} 
E_3=\sum_{s=1}^{J}\bar a_{s,k}(r)e_s, \quad 
E_4=\left\{\begin{array}{ll} 
O(M_k^{-\frac{18-\epsilon}{n-2}}r^{9-\frac{\epsilon}2-n}),\quad n=10,\\ \\ 
O(M_k^{-\frac{20-\epsilon}{n-2}}r^{10-\frac{\epsilon}2-n}),\quad n\ge 11. 
\end{array} 
\right. 
\label{dec9e9} 
\end{equation} 
where $e_s=e_s(\theta)$, independent of $k$, is a homogeneous 
spherical harmonic of degree $s$, $J$ is a positive integer, 
 $\bar a_{s,k}$ 
satisfies 
$$ 
|\bar a_{s,k}(r)| 
=O(M_k^{-\frac{16-\epsilon}{n-2}}r^{3-n})+O(M_k^{-\frac{16}{n-2}}r^{4-n}). 
$$ 
From now on we say a term is $E_3$ or $E_4$ if it is of the form in 
(\ref{dec9e9}). The main result in this subsection is 
\begin{prop} 
\label{propnew1} 
 For $n\ge 10$, $\frac 12\le \lambda\le 2$, 
$$ E_{\lambda}=E_1+E_2+E_3+E_4,\qquad \lambda\le |y|\le \frac 12 
M_k^{ \frac 2{n-2} }. 
$$ 
\end{prop} 
 
\noindent{\bf Proof.}\ 
First by (\ref{nov19e6}) we have 
\begin{equation} 
\bigg|\nabla^l \left[v_k^{\lambda}- 
(U^{\lambda}+F^{(3)}_{\lambda})\right]\bigg|\le C 
M_k^{-\frac{12}{n-2}}|y|^{2-n},\quad l=0,1,2, 
\label{dec9e2} 
\end{equation} 
where $F^{(3)}_{\lambda}$ is the Kelvin transformation of $F^{(3)}$. 
Note that (\ref{dec9e2}) holds over the whole $\Sigma_{\lambda}$. 
Similarly we can define $F^{(2)}_{\lambda}$ as the Kelvin 
transformation of $F^{(2)}$ and we shall use this expression: 
$$F^{(2)}_{\lambda}= 
-c(n)\widetilde R^{(2)}(\theta) 
f^{\lambda}_2(r)M_k^{-\frac 8{n-2}}. 
$$ 
where 
$$ 
\widetilde R^{(2)}(\theta)= 
\sum_{i<j}\partial_{ij}R\theta_i\theta_j+\sum_i\frac{\partial_{ii}R}{2} 
(\theta_i^2-\frac 1n) 
=\sum_{i<j}\partial_{ij}R\theta_i\theta_j+ 
\frac 12 \sum_i (\partial_{ii}R) 
\theta_i^2+O(|W|^2). 
$$ 
and $f^{\lambda}_{2}(r)=(\frac{\lambda}r)^{n-2}f_{2}(\lambda^2/r)$ 
is the Kelvin transformation of $f_{2}$. $f_3^{\lambda}$ is 
understood similarly. Since $0\le f_2(r),f_3(r)\le Cr$ for $0\le 
r\le 1$, we have 
\begin{equation} 
|f_2^\lambda(r)|+|f_3^\lambda(r)|\le Cr^{1-n},\qquad 
F^{(3)}_{\lambda}=O(M_k^{-\frac{8}{n-2}})r^{1-n}. 
\label{926e1} 
\end{equation} 
 
Now we consider $(\bar b_i\partial_i+\bar d_{ij}\partial_{ij})v_k^{\lambda}$. 
Since $U^{\lambda}$ is radially symmetric, 
we have, using (\ref{925e1}), 
(\ref{dec9e1}), 
(\ref{dec9e2}), 
\begin{eqnarray*} 
&& (\bar b_i\partial_i+\bar d_{ij}\partial_{ij})v_k^{\lambda} =(\bar 
b_i\partial_i+\bar d_{ij}\partial_{ij}) 
F^{(3)}_{\lambda} +E_4\nonumber\\ 
&=&-c(n) 
\bigg\{(\bar b_i\partial_i+\bar d_{ij}\partial_{ij}) 
\widetilde R^{(2)}(\theta) 
\bigg\} 
f^{\lambda}_{2}(r)M_k^{-\frac{8}{n-2}}\nonumber \\ 
&&- 
c(n)(\bar b_i\partial_i+\bar d_{ij}\partial_{ij})\sum_{p=1}^ 
{l_3}\tilde R^{(3)}_{p}e_{3p}f_3^{\lambda} 
M_k^{-\frac{10}{n-2}}+E_4. 
\end{eqnarray*} 
 
For any smooth function $a(\theta)$, 
\begin{equation} 
\int_{S^{n-1}}(\bar b_i\partial_i+\bar d_{ij}\partial_{ij}) 
a(\theta)= 
\int_{S^{n-1}}(\Delta_{g_k}-\Delta) a(\theta)=0. 
\label{P2} 
\end{equation} 
Expanding $\bar b_i(y)$ and $\bar d_{ij}(y)$ to the fourth and the 
fifth order respectively and using (\ref{8-2}), we have 
$$ 
(\bar b_i\partial_i+\bar d_{ij}\partial_{ij}) 
\widetilde R^{(2)}(\theta)= \sum_{l=1}^7 a_{l,k}(r) 
P_l(\theta)+O(M_k^{ -\frac{12}{n-2} }r^4), 
$$ 
$$ 
(\bar b_i\partial_i+\bar d_{ij}\partial_{ij})\sum_{p=1}^{l_3}\widetilde 
R^{(3)}_{p}e_{3p} 
=\sum_{l=1}^5 b_{l,k}(r)P_l(\theta)+O(M_k^{ -\frac{12}{n-2} }r^4), 
$$ 
where $a_{l,k}(r)$ and $b_{l,k}(r)$ are radial functions satisfying 
\begin{equation} 
|a_{l,k}(r)|+|b_{l,k}(r)| 
=O(M_k^{-\frac{8-\epsilon}{n-2}}r)+O(M_k^{-\frac{8}{n-2}}r^2), 
\label{P1} 
\end{equation} 
while $P_l(\theta)$ is a homogeneous polynomial in $\theta$ of 
degree $l$ and is also independent of $k$. Consequently, using 
(\ref{P2}), we have 
$$ 
(\bar b_i\partial_i+\bar d_{ij}\partial_{ij}) 
\widetilde R^{(2)}(\theta) 
= \sum_{l=1}^7 a_{l,k}(r)e_l(\theta)+O(M_k^{ -\frac{12}{n-2} }r^4), 
$$ 
$$ 
(\bar b_i\partial_i+\bar d_{ij}\partial_{ij})\sum_{p=1}^{l_3}\widetilde 
R^{(3)}_{p}e_{3p} 
=\sum_{l=1}^7 b_{l,k}(r)e_l(\theta)+O(M_k^{ -\frac{12}{n-2} }r^4), 
$$ 
where $e_l(\theta)$, independent of $k$, is 
a homogeneous spherical harmonic of degree $l$, and 
$ 
a_{l,k}(r)$ and $b_{l,k}(r)$, independent of $\theta$, satisfy (\ref{P1}). 
Consequently, using also  (\ref{926e1}), 
$$ 
 (\bar b_i\partial_i+\bar d_{ij}\partial_{ij})v_k^{\lambda}= 
(\bar b_i\partial_i+\bar d_{ij}\partial_{ij})F^{(3)}_{\lambda} 
+E_4=E_3+E_4. 
$$ 
 
 Similarly we can show that 
$$(\frac{\lambda}{|y|})^{n+2}(\bar b_i(y^{\lambda})\partial_i 
v_k(y^{\lambda})+\bar d_{ij}(y^{\lambda})\partial_{ij}v_k(y^{\lambda})) 
=E_3+E_4.$$ 
 
We have discussed the minor terms in $E_{\lambda}$, by (\ref{eQ}), 
the main term in $E_{\lambda}$ is 
$$\displaystyle{\bar 
c(y)v_k^{\lambda}-(\frac{\lambda}{|y|})^{n+2}\bar c(y^{\lambda}) 
v_k(y^{\lambda})}.$$ 
 
We shall use the following two expansions of $\bar c$ according to 
circumstances.

First we know 
\begin{eqnarray} 
\bar c&=&c(n)\sum_{l=2}^{7}\tilde R^{(l)}(\theta)r^l 
M_k^{-\frac{4+2l}{n-2}}+c(n)\sum_{s=0}^{2}\bar R^{(2s+2)}r^{2s+2} 
M_k^{-\frac{8+4s}{n-2}}\nonumber \\ 
&&\qquad +O(M_k^{-\frac{20}{n-2}}r^8). \qquad \qquad \lambda<r<\frac 
12 M_k^{ \frac 2{n-2}}. \label{comega2} 
\end{eqnarray} 
 
On the other hand, using (\ref{SS5}), (\ref{7-1new}) 
 and the rate of $|W|$, 
\begin{equation} 
\bar c=c^{(3)}+O(M_k^{-\frac{12}{n-2}})r^4 
+O(M_k^{ -\frac{14-\epsilon}{n-2} }r^3), \qquad \qquad 
 \lambda<r<\frac 12 M_k^{ \frac 2{n-2}}. 
\label{comega3} 
\end{equation} 
where 
\begin{equation} 
c^{(3)}= 
c^{(2)}+ 
c(n)\sum_{p=1}^{l_3}\tilde R^{(3)}_{p}e_{3p}r^3M_k^{-\frac{10}{n-2}}, 
\label{comega1} 
\end{equation} 
$$ 
c^{(2)}=c(n) \bigg (\sum_{i<j}\partial_{ij}R\theta_i\theta_j+\frac 12 
\sum_i\partial_{ii}R\theta_i^2\bigg ) 
r^2M_k^{-\frac{8}{n-2}}. 
$$ 
 
Using (\ref{dec9e2}), we have 
\begin{eqnarray} 
&&\bar cv_k^{\lambda}-(\frac{\lambda}r)^{n+2}\bar c(y^{\lambda}) 
v_k(y^{\lambda}) \nonumber \\ 
&=&\bar cU^{\lambda}-(\frac{\lambda}r)^{n+2}\bar c(y^{\lambda})U(y^{\lambda}) 
+c^{(3)}F^{(3)}_{\lambda}-(\frac{\lambda}{r})^{n+2} c^{(3)} 
(y^{\lambda}) 
F^{(3)}(y^{\lambda})+E_4. 
\label{926e2} 
\end{eqnarray} 
 
To estimate the first two terms 
we use the  expression  (\ref{comega2}) 
 of $\bar c$. 
 
\begin{eqnarray*} 
&&\bar cU^{\lambda}-(\frac{\lambda}r)^{n+2}\bar c(y^{\lambda})U(y^{\lambda}) \nonumber \\ 
&=&c(n)U^{\lambda} 
\bigg (\sum_{l=2}^{7}\tilde R^{(l)}r^lM_k^{-\frac{4+2l}{n-2}}+ 
\sum_{s=0}^{2}\bar R^{(2s+2)}M_k^{-\frac{8+4s}{n-2}}r^{2s+2}\bigg ) \nonumber \\ 
&&-c(n)(\frac{\lambda}r)^{n+2}U(y^{\lambda})\bigg ( 
\sum_{l=2}^{7}\tilde R^{(l)}(\frac{\lambda^2}{r})^l 
M_k^{-\frac{4+2l}{n-2}}+\sum_{s=0}^{2}\bar R^{(2s+2)} 
M_k^{-\frac{8+4s}{n-2}}(\frac{\lambda^2}{r})^{2s+2}\bigg ) \nonumber \\ 
&&+E_4 \nonumber \\ 
&=&c(n)U^{\lambda} 
\sum_{s=0}^{2}\bar R^{(2s+2)}M_k^{-\frac{8+4s}{n-2}}r^{2s+2} 
(1-(\frac{\lambda}{r})^{4s+8}) 
+E_2 +E_4. 
\end{eqnarray*} 
 
For the third and the fourth terms of (\ref{926e2}) we use (\ref{comega3}), 
(\ref{comega1}) and (\ref{926e1}) to obtain 
\begin{eqnarray*} 
&& c^{(3)}F^{(3)}_{\lambda}-(\frac{\lambda}r)^{n+2} 
 c^{(3)}(y^{\lambda})F^{(3)}(y^{\lambda})\\ 
&=& 
\bigg (c^{(3)}-(\frac{\lambda}{r})^4 
c^{(3)}(y^{\lambda})\bigg )F^{(3)}_{\lambda}\\ 
&=& 
\bigg (c^{(2)}-(\frac{\lambda}{r})^4 
c^{(2)}(y^{\lambda})\bigg )F^{(2)}_{\lambda}+E_3\\ 
&=&-c(n)^2 \bigg (\sum_{i<j}\partial_{ij}R\theta_i\theta_j+\frac 12 
\sum_i\partial_{ii}R\theta_i^2\bigg )^2 
r^2(1-(\frac{\lambda}{r})^{8}) 
f^{\lambda}_{2} 
M_k^{-\frac{16}{n-2}}\\ 
&& +E_3+E_4. 
\end{eqnarray*} 
 
In the expansion of 
 the product in the last equality, we use the fact that 
homogeneous polynomials of 
 degree $2$ and $3$ are orthogonal to each other --- when a term 
has average 0 on $\Bbb S^{n-1}$ it contributes a term $E_3$.

The term that needs to be evaluated is $ \bigg 
(\sum_{i<j}\partial_{ij}R\theta_i\theta_j+\frac 12 
\sum_i\partial_{ii}R\theta_i^2\bigg )^2. $ It is elementary to 
verify the following identities: 
$$\frac{1}{|S^{n-1}|}\int_{S^{n-1}}\theta_i^2\theta_j^2 
=\frac{1}{n(n+2)}\quad i\neq j,\qquad 
\frac{1}{|S^{n-1}|}\int_{S^{n-1}} (\theta_i)^4= \frac 3{  n(n+2) }. 
$$

In the following we often write a polynomial $P(\theta)$  of 
degree less or equal to $7$ as the sum of 
$\frac 1{|S^{n-1}|}\int_{S^{n-1}} P(\theta)$ 
and $\sum_{p=1}^7 C_p e_p(\theta)$ where 
$e_p(\theta)$ are homogeneous spherical harmonics 
of degree $p$. 
 
By the above, we write 
\begin{eqnarray*} 
&& \left(  \sum_{i<j}\partial _{ij}R \theta_i \theta_j 
+\frac 12 \sum_i \partial _{ii}R \theta_i^2\right)^2 
\\ 
&=& \sum_{i<j} (\partial _{ij}R)^2 \theta_i ^2 \theta_j^2 
+\frac 14 \left( \sum_i \partial _{ii}R \theta_i^2\right)^2 
+\sum_{p=1}^7 C_p e_p(\theta)\\ 
&=&  \sum_{i<j} (\partial _{ij}R)^2 \theta_i ^2 \theta_j^2 
+\frac 14 \sum_{i\ne j} \partial_{ii}R 
\partial_{jj}R \theta_i^2 \theta_j^2+ 
\frac 14 \sum_{i}(\partial _{ii}R)^2 \theta_i^4 
+\sum_{p=1}^7 C_p e_p(\theta)\\ 
&=& \frac 1{ n(n+2) }  \sum_{i<j} (\partial _{ij}R)^2 
+\frac 1{ 4n(n+2) }  \sum_{i\ne j} \partial_{ii}R 
\partial_{jj}R +\frac 
 3{ 4n(n+2) }    \sum_{i}(\partial _{ii}R)^2 
+\sum_{p=1}^7 C_p e_p. 
\end{eqnarray*} 
Note that 
$$\sum_{i\neq j}\partial_{ii}R\partial_{jj}R+\sum_{i=1}^n(\partial_{ii}R)^2= 
(\sum_{i}\partial_{ii}R)^2=O(|W|^4)=O(M_k^{-\frac{16-\epsilon}{n-2}}).$$ 
We have 
\begin{eqnarray} 
 &&\left(  \sum_{i<j}\partial _{ij}R \theta_i \theta_j 
+\frac 12 \sum_i \partial _{ii}R \theta_i^2\right)^2\nonumber\\ 
&=&\frac 1{ 2n(n+2) } 
[\sum_{i<j}2(\partial_{ij}R)^2+\sum_{i}(\partial_{ii}R)^2 ] 
+O(M_k^{-\frac{16-\epsilon}{n-2}}) +\sum_{p=1}^7C_pe_p. 
\label{14-1new} 
\end{eqnarray} 
 
Thus 
\begin{eqnarray*} 
&& c^{(3)}F^{(3)}_{\lambda}-(\frac{\lambda}r)^{n+2} 
 c^{(3)}(y^{\lambda})F^{(3)}(y^{\lambda})\\ 
&=& 
-\frac{c(n)^2}{2n(n+2)}[\sum_{i<j}2(\partial_{ij}R)^2+\sum_{i}(\partial_{ii}R)^2 
]r^2(1-(\frac{\lambda}{r})^8)f_2^{\lambda}M_k^{-\frac{16}{n-2}} 
+E_3+E_4. 
\end{eqnarray*} 
Proposition \ref{propnew1} follows from the above. $\Box$

\subsection{Construction of auxiliary  functions 
and proof of (\ref{eq4})} The goal of this subsection is to finish 
the proof of (\ref{eq4}) by finding a contradiction to (\ref{5-0}). 
Before we construct the auxiliary 
 functions, we discuss two relatively minor terms in 
$E_1$. Recall that $\bar R^{(2)}=-\frac{|W|^2}{12n}.$ The following 
properties of conformal normal coordinates are established in 
\cite{HV}: If $W=0$, then 
 $ 
  R_{abcd}=0$ and, for some constant 
$c_1(n)>0$, 
 $\bar R^{(4)}= 
 -c_1(n)|R_{abcd,e}|^2$; 
if $W=0$ and $\nabla W=0$, then 
$R_{abcd,e}=0$. 
Examining the proofs there, 
we arrive at 
$$ 
|R_{abcd}|=O(1)|W|,\quad |R_{abcd,e}|=O(|W|)+O(|\nabla_g W|), 
$$ 
and 
$$\bar R^{(4)}=-c_1(n) 
|\nabla R_{abcd}|^2+O(|W|)O(|\nabla R_{abcd}|)+O(|W|^2).$$Thus 
\begin{equation} 
\bar R^{(4)}\le -\frac 12 c_1(n)|\nabla R_{abcd}|^2+O(|W|^2)= 
-\frac 12 c_1(n)|\nabla R_{abcd}|^2+O(M_k^{-\frac{8-2\epsilon}{n-2}}). 
\label{dec10e5} 
\end{equation} 
So 
\begin{eqnarray} 
E_1&\le& c(n)U^{\lambda}\bar R^{(6)} r^6 M_k^{ -\frac{16}{n-2} } ( 
1-(\frac{\lambda}{r})^{16}) 
\nonumber\\ 
&&-\frac{c(n)^2}{2n(n+2)}\bigg (\sum_{i<j} 
2(\partial_{ij}R)^2+\sum_{i}(\partial_{ii}R)^2\bigg )(1-(\frac{\lambda}{r})^{8})r^2f_2^{\lambda} 
M_k^{-\frac{16}{n-2}}. 
\label{15-6new} 
\end{eqnarray}

Recall that the most important requirement for the test function 
$h_{\lambda}$ is (\ref{20-1new}), for which we  construct 
$h_{\lambda}$ as the sum of four test functions $h_1,...,h_4$. Each of 
the first three functions is constructed with respect to $\Delta$ 
and $V_{\lambda}$ (a radial function, to be defined later), rather 
than $\Delta_{g_k}-\bar c$ and $\xi$. So even though each of them 
cancels a major part of $E_{\lambda}$, they also create some minor 
extra errors because of the difference between $\Delta, V_{\lambda}$ 
and $\Delta_{g_k}-\bar c,\xi$. Eventually all these new error terms 
will be put together and be controlled by $h_4$. 
 
For the convenience of our discussion, we define 
$$\bar l=6\,\, \mbox{if}\quad n=10,\quad \bar l=7 \,\, 
 \mbox{if}\quad n= 
 11.$$ 
 
To cancel the 
term $\displaystyle{c(n)U^{\lambda}\sum_{l=2}^{\bar l}\tilde R^{(l)} 
M_k^{-\frac{4+2l}{n-2}}r^l(1-(\frac{\lambda}{r})^{2l+4})}$ in $E_{\lambda}$ we 
use (\ref{15-4new}). 
  The 
$f_{2,\lambda}$ defined in Appendix C is to deal with $\tilde 
R^{(2)}$. 
 
Let 
$$V_{\lambda}(r):=n(n+2)\int_0^1(tU+(1-t)U^{\lambda})^{\frac 
4{n-2}}dt.$$ 
 
For $3\le l\le \bar l$, 
let $\lambda_p>0$ be the eigenvalue corresponding to $e_p$, we consider 
\begin{equation} 
\left\{\begin{array}{ll} 
f_{p\lambda l}''(r)+\frac{n-1}rf_{p\lambda l}'(r)+ 
(V_{\lambda}-\frac{\lambda_p}{r^2})f_{p\lambda l}(r)=-r^lU^{\lambda}(r) 
(1-(\frac{\lambda}{r})^{2l+4}),\\ 
\qquad \qquad \qquad \qquad \qquad \lambda<r<2M_k^{\frac{2}{n-2}}, 
\quad 3\le l\le \bar l.\\ 
f_{p\lambda l}(\lambda)=f_{p\lambda l}(2M_k^{\frac{2}{n-2}})=0. 
\end{array} 
\right. 
\label{15-1new} 
\end{equation}

By 
Proposition 6.1 in Appendix A of \cite{LiZhang1}, 
there exists some small $\epsilon_4=\epsilon_4(n)>0$ such that 
for $\lambda \in [1-\epsilon_4, 1+\epsilon_4]$, 
equation (\ref{15-1new}) has a unique classical solution 
satisfying 
$$0\le f_{p\lambda l}(r)\le Cr^{l+4-n},\quad 3\le l\le \bar l,\quad \lambda\le r\le 
2M_k^{\frac{2}{n-2}}.$$

Let 
$$h_1=c(n) 
\widetilde R^{(2)}(\theta) 
f_{2,\lambda}M_k^{-\frac 8{n-2}}.$$ 
This is a major part of $h_{\lambda}$. Let 
$$h_2=c(n)\sum_{l=3}^{\bar l}\sum_{p=1}^{I_l}\tilde R_{lp}e_p(\theta)f_{p\lambda l}(r) 
M_k^{-\frac{4+2l}{n-2}},$$ where $\tilde R_{lp}$ and $e_p(\theta)$ 
are the ones in (\ref{15-4new}). By the definitions of $h_1$ and 
$h_2$, we have 
$$\Delta (h_1+h_2)+V_{\lambda}(h_1+h_2)=-E_2.$$

The extra error terms created by $h_1$ are 
$$\bigg (\bar b_i\partial_i+\bar d_{ij}\partial_{ij}-\bar c+(n(n+2)\xi^{\frac{4}{n-2}} 
-V_{\lambda})\bigg )h_1.$$ We need to estimate the above in 
$O_{\lambda}$, note that by definition $h_1, h_2=\circ(1)|y|^{2-n}$ 
in $\Sigma_{\lambda}$. For $(\bar b_i\partial_i+\bar 
d_{ij}\partial_{ij})h_1$, we just analyze it the same
way as analyzing $(\bar 
b_i\partial_i+\bar d_{ij}\partial_{ij})F^{(3)}_{\lambda}$ to obtain
$$ 
(\bar b_i\partial_i+\bar d_{ij}\partial_{ij})h_1= 
\widetilde E_3+E_4. 
$$ 
Here and in the following, $\widetilde E_3$ denotes 
a term of the form 
$$ 
\widetilde E_3=\sum_{s=1}^J\bar c_{s,k}(r)e_s 
$$ 
with  $\bar c_{s,k}(r)$ depending only on $r$ and satisfying 
\begin{equation} 
\bar c_{s,k}(r)=O(M_k^{-\frac{16-\epsilon}{n-2}}r^{7-n}) 
+O(M_k^{-\frac{16}{n-2}}r^{8-n}). 
\label{dec10e11} 
\end{equation}

Next we consider $-\bar ch_1$. By the definition of $h_1$,
 (\ref{comega3}) 
and (\ref{comega1}) we have 
$$-\bar ch_1=- c^{(3)}h_1+O(M_k^{-\frac {20}{n-2}})r^{10-n}.$$

The major part to contribute to $E_1$ is, 
using (\ref{14-1new}), 
\begin{eqnarray*} 
- c^{(2)}h_1 
&=&-c(n)^2\bigg (\sum_{i<j}\partial_{ij}R\theta_i\theta_j+\frac 12\sum_i\partial_{ii}R 
\theta_i^2\bigg )^2 
r^2f_{2,\lambda}M_k^{-\frac{16}{n-2}}
+E_4\nonumber \\ 
&=&-\frac{c(n)^2}{2n(n+2)} \bigg 
(\sum_{i<j}2(\partial_{ij}R)^2+\sum_i(\partial_{ii}R)^2 \bigg 
)r^2f_{2,\lambda}M_k^{-\frac{16}{n-2}} +\widetilde E_3+E_4. 
\end{eqnarray*}

Since $\widetilde R^{(2)}(\theta)$ is orthogonal to $e_{3p}(\theta)$ 
\begin{equation} 
(c^{(3)}-c^{(2)})h_1=\widetilde E_3. 
\label{17-1new} 
\end{equation} 
 
For $-\bar ch_1+(n(n+2)\xi^{\frac 4{n-2}}-V_{\lambda})h_1$ we use 
(\ref{nov19e6}) over $|y|\le M_k^{\frac{16-\epsilon}{(n-2)^2}}$ and 
$$n(n+2)\xi^{\frac 4{n-2}}-V_{\lambda}=O(r^{-4}),\quad r\in 
(M_k^{\frac{16-\epsilon}{(n-2)^2}},  \frac 1{\sqrt k}M_k^{\frac 
2{n-2}}).$$ 
 
First for $r\le M_k^{\frac{16-\epsilon}{(n-2)^2}}$, 
using (\ref{nov19e6}) and (\ref{dec9e2}), 
\begin{eqnarray} 
&&n(n+2)\xi^{\frac 4{n-2}}=n(n+2)\int_0^1(t(U+a)+(1-t)(U^{\lambda}+b)) 
^{\frac 4{n-2}}dt \nonumber \\ 
&=&V_{\lambda}+\frac{4n(n+2)}{n-2} 
\int_0^1(tF^{(3)}+(1-t)F^{(3)}_{\lambda})(tU+(1-t)U^{\lambda})^{\frac{6-n}{n-2}}dt 
\nonumber \\ 
&&+O(M_k^{-\frac {12}{n-2}}r^{2+\bar a}).\quad \lambda\le r\le 
M_k^{\frac{16-\epsilon}{(n-2)^2}}, 
\label{dec10e15} 
\end{eqnarray} 
where $a 
:= v_k-U=F^{(3)}+O(M_k^{-\frac {12}{n-2}}r^{8-n+\bar a})$ and 
$b:= v_k^\lambda-U^\lambda 
=F^{(3)}_{\lambda}+O(M_k^{-\frac {12}{n-2}}r^{2-n})$. 
 
Since $\widetilde R^{(2)}(\theta)$ is orthogonal to $e_{3p}(\theta)$ 
 
\begin{eqnarray*} 
&&\bigg (\int_0^1(tF^{(3)}+(1-t)F^{(3)}_{\lambda}-tF^{(2)}-(1-t)F^{(2)}_{\lambda}) 
(tU+(1-t)U^{\lambda})^{\frac{6-n}{n-2}}dt\bigg )h_1\\ 
&=&\widetilde E_3. 
\end{eqnarray*}

Clearly, 
$$ 
(tF^{(2)}+(1-t)F^{(2)}_{\lambda})(tU+(1-t)U^{\lambda})^{\frac{6-n}{n-2}}h_1 
\le 0. 
$$

So 
$$ 
\left(n(n+2)\xi^{\frac 4{n-2}}-V_{\lambda}\right) 
h_1\le 
\widetilde E_3+E_4, 
\qquad \quad \lambda\le r\le 
M_k^{\frac{16-\epsilon}{(n-2)^2}}. 
$$

We use 
$$|(n(n+2)\xi^{\frac{4}{n-2}}-V_{\lambda})h_1| 
\le O(M_k^{-\frac{8}{n-2}}r^{2-n}) =E_4,\qquad 
M_k^{\frac{16-\epsilon}{(n-2)^2}}\le r\le  \frac 1{\sqrt k} 
M_k^{\frac{2}{n-2}}. 
$$

Now we estimate 
$$(\bar b_i\partial_i+\bar d_{ij}\partial_{ij}-\bar c 
+n(n+2)\xi^{\frac{4}{n-2}}-V_{\lambda})h_2.$$ 
As usual, 
$$(\bar b_i\partial_i+\bar d_{ij}\partial_{ij})h_2= 
\widetilde E_3, 
$$ 
$$ 
-\bar ch_2=- c^{(2)}h_2+ 
E_4 
=- c^{(2)} 
c(n)\left( 
\sum_{p=1}^{I_l}\tilde R_{lp}e_p(\theta)f_{p\lambda l}(r) 
 M_k^{-\frac{4+2l}{n-2}}\right)\bigg|_{l=3}+ 
E_4. 
$$ 
Compare (\ref{7-1new}) and  (\ref{15-4new}), 
we deduce from above, using the orthogonality of 
$\widetilde R^{(2)}$ and $e_{3p}(\theta)$ and the decay of 
$W$, 
\begin{equation} 
-\bar ch_2=\widetilde E_3+E_4. 
\label{18-1new} 
\end{equation}

By (\ref{dec10e15}) 
\begin{eqnarray*} 
&&n(n+2)\xi^{\frac{4}{n-2}}-V_{\lambda}\\ 
&=&\frac{4n(n+2)}{n-2} 
\int_0^1(tF^{(2)}+(1-t)F^{(2)}_{\lambda})(tU+(1-t)U^{\lambda})^{\frac{6-n}{n-2}}dt\\ 
&&+O(M_k^{-\frac{10}{n-2}}r)+O(M_k^{-\frac{12}{n-2}}r^{2+\bar a}),\quad 
\lambda<r<M_k^{\frac{16-\epsilon}{(n-2)^2}}. 
\end{eqnarray*} 
Therefore, as in the derivation of 
(\ref{18-1new}),

\begin{eqnarray*} 
&&(n(n+2)\xi^{\frac{4}{n-2}}-V_{\lambda})h_2\\ 
&=&\frac{4n(n+2)}{n-2} 
h_2 \int_0^1(tF^{(2)}+(1-t)F^{(2)}_{\lambda}) 
(tU+(1-t)U^{\lambda})^{\frac{6-n}{n-2}}dt+E_4\\ 
&=& \frac{4n(n+2)}{n-2} 
\left( 
c(n) 
\sum_{p=1}^{I_l}\tilde R_{lp}e_p(\theta)f_{p\lambda l}(r) 
 M_k^{-\frac{4+2l}{n-2}}\right)\bigg|_{l=3}\\ 
&& \cdot 
\int_0^1(tF^{(2)}+(1-t)F^{(2)}_{\lambda}) 
(tU+(1-t)U^{\lambda})^{\frac{6-n}{n-2}}dt+E_4\\&=& \widetilde E_3+E_4 
\qquad\qquad\qquad\lambda<r<M_k^{\frac{16-\epsilon}{(n-2)^2}}. 
\end{eqnarray*} 
The first term on the right is of the form $\sum_{s=1}^J\bar 
c_{s,k}(r)e_s$ and $\bar 
c_{s,k}(r)=O(M_k^{-\frac{18}{n-2}}r^{7-n})$. So after this term is 
extended to $r\le \frac 1{\sqrt k} M_k^{\frac{2}{n-2}}$, it can be 
combined with $E_3$. The extended part has a good decay. 
 
For $M_k^{\frac{16-\epsilon}{(n-2)^2}}\le r\le \frac 1{\sqrt k} 
M_k^{\frac{2}{n-2}}$, we have 
$$|(n(n+2)\xi^{\frac{4}{n-2}}-V_{\lambda})h_2| 
\le O(M_k^{-\frac{10}{n-2}}r^{3-n})=E_4. 
$$ 
Recall that our purpose is to obtain (\ref{20-1new}). By putting 
$h_1$ and $h_2$ together and using (\ref{15-6new}), we have 
\begin{eqnarray} 
&&(\Delta+\bar b_i\partial_i+\bar d_{ij}\partial_{ij} 
-\bar c+n(n+2)\xi^{\frac{4}{n-2}})(h_1+h_2)+E_1+E_2\nonumber \\ 
&\le  & c(n)U^{\lambda}\bar R^{(6)}M_k^{-\frac{16}{n-2}}r^{6} 
(1-(\frac{\lambda}{r})^{16})\nonumber \\ 
&-&\frac{c(n)^2}{2n(n+2)}\sum_{i,j}(\partial_{ij}R)^2 
r^2 \bigg ((1-(\frac{\lambda}{r})^{8})f_2^{\lambda}+f_{2,\lambda}\bigg ) 
M_k^{-\frac {16}{n-2}}+\widetilde E_3+E_4. 
\label{dec10e18} 
\end{eqnarray} 
 
For $\bar R^{(6)}$ we use (\ref{dec10e8}), we also know the lower 
bounds for $f_2$ and $f_{2,\lambda}$, respectively (see 
(\ref{nov19e1}) and (\ref{mar10e2})). These three estimates are in 
the appendix and are sufficient for $h_1+h_2$ to cancel the major 
part of $E_{\lambda}$. In fact, first by (\ref{nov19e1}) and 
(\ref{mar10e2}) 
\begin{eqnarray*} 
&&\frac{c(n)}{2n(n+2)} \bigg 
((1-(\frac{\lambda}{r})^{8})f_2^{\lambda}+f_{2,\lambda}\bigg )\\ 
&\ge&U^{\lambda}r^4(1-(\frac{\lambda}{r})^{16})\frac{1}{8(n+4)(n+2)n} 
(\frac{n-8}{n-2}-\frac{49}{20n^2} +\epsilon) 
\end{eqnarray*} 
where we have used the following inequality that holds only for 
$n=10,11$. 
$$ 
\frac{1}{8(n+4)(n+2)n}(\frac{n-8}{n-2}-\frac{49}{20n^2} 
+\epsilon) 
\le \frac{c(n)}{2n(n+2)}\frac{1}{6(n-4)}, 
$$ 
Thus, by using (\ref{dec10e8}), we deduce from (\ref{dec10e18}) that 
\begin{eqnarray*} 
&&(\Delta+\bar b_i\partial_i+\bar d_{ij}\partial_{ij} 
-\bar c+n(n+2)\xi^{\frac{4}{n-2}})(h_1+h_2)\nonumber \\ 
&\le &-E_1-E_2+\widetilde E_3+E_4. 
\quad \mbox{in}\quad \Sigma_{\lambda}. 
\end{eqnarray*} 
 
Next we construct test functions to control $\widetilde E_3$ and 
$E_4$. The $h_3$ to be constructed later will create much minor 
error terms than before. Then eventually all the minors terms will 
be controlled by $h_4$. Let 
 $\tilde f_{s \lambda}$ be the solution of 
$$\left\{\begin{array}{ll} 
\tilde f_{s\lambda}''(r)+\frac{n-1}r\tilde f_{s\lambda}'(r) 
+(V_{\lambda}- 
\frac{\lambda_s}{r^2})\tilde f_{s\lambda}(r)=-\bar c_{s,k}(r), 
\quad \lambda<r<M_k^{\frac{2}{n-2}},\\  \\ 
\tilde f_{s\lambda}(\lambda)=\tilde 
f_{s\lambda}(M_k^{\frac{2}{n-2}})=0, 
\end{array} 
\right. 
$$ 
and let 
$$h_3:=\sum_{s=1}^{J}\tilde f_{s\lambda}(r)e_s.$$ 
 
By (\ref{dec10e11}), $\bar 
c_{s,k}(r)=O(M_k^{-\frac{14}{n-2}}r^{7-n}).$ Consequently $|\tilde 
f_{s\lambda}(r)|\le CM_k^{-\frac{14}{n-2}}r^{9-n}.$ Therefore 
$$(\Delta +V_{\lambda})h_3=-\widetilde E_3.$$ 
By the estimates of $\bar b_i$, $\bar d_{ij}$ and $h_3$, etc we obtain 
$$(\bar b_i\partial_i+\bar d_{ij}\partial_{ij}-\bar c 
+n(n+2)\xi^{\frac{4}{n-2}}-V_{\lambda})h_3= E_4.$$ Finally we 
define, for $Q>>1$ that 
$$h_4(r)=\left\{\begin{array}{ll} 
QM_k^{-\frac{18-2\sqrt{\epsilon}}{n-2}}f_{n,n-9+\sqrt{\epsilon}} 
(\frac r{\lambda}),&\quad n=10,\\ 
\\ 
QM_k^{-\frac{20-2\sqrt{\epsilon}}{n-2}}f_{n,n-10+\sqrt{\epsilon}} 
(\frac r{\lambda}),&\quad n=11. 
\end{array} 
\right. 
$$ 
where $f_{n,\alpha}$ is defined in \cite{LiZhang1}. Let 
$h_\lambda:=h_1+h_2+h_3+h_4$, then $(\ref{20-1new})$ is obtained. 
This $h_{\lambda}$ satisfies all the requirements for the test 
function to make the method of moving spheres work. Then the 
standard moving sphere argument leads to the following 
conclusion: 
$$\min_{|y|\le r}v_k\le (1+\epsilon)U(r),\quad 0\le r\le \frac 
1{\sqrt k}M_k^{\frac 2{n-2}}$$ where $\epsilon$ is an arbitrary 
small positive constant. Then following the argument in 
\cite{LiZhang1} one gets a contradiction to (\ref{5-0}). (\ref{eq4}) 
is established. 
 
\subsection{Vanishing rates of the Weyl tensor and 
the completion of the proof of Theorem \ref{thm1} } 
 
In this subsection we use (\ref{eq4}) to prove 
(\ref{W1}) and (\ref{W2}), 
the vanishing rates  of 
the  Weyl tensor and its covariant derivatives at the blow up point. 
By (\ref{eq4}), 
\begin{equation} 
v_k(y)\le C U(y), \qquad 
|y|\le \delta M_k^{ \frac 2{n-2} }. 
\label{21-new9} 
\end{equation} 
This estimate leads to an improved estimate
of $v_k$ than that in 
Proposition \ref{propnew1}. 
\begin{prop} There exists $\delta'>0$, independent of $k$, such that 
\begin{eqnarray*} 
&&\bigg|\nabla^l\left(
v_k-( v^{(1)}+ M_k^{ -\frac 8{n-2} }
v^{(2)}+ M_k ^{ -\frac {10}{n-2} }
v^{(3)})\right)\bigg|
\nonumber\\
&=&O(M_k^{-\frac{12}{n-2}})(1+|y|)^{8-n-l}, 
\qquad\qquad \qquad |y| \le \delta' M_k^{\frac{2}{n-2}}, 
\ l=0,1,2. 
\end{eqnarray*} 
\label{E1-1} 
\end{prop} 
 
\noindent {\bf Proof.}\ 
Write 
$$
E^{(3)}:= v_k-(U+F^{(3)}).
$$
We only need to prove that
\begin{equation}
|\nabla^l E^{(3)}|=O(M_k^{-\frac{12}{n-2}})(1+|y|)^{8-n-l},
\qquad |y| \le \delta' M_k^{\frac{2}{n-2}},
\ l=0,1,2.
\label{2omega82}
\end{equation}

It follows from (\ref{ab1}) and (\ref{dec8e1}) that 
\begin{equation} 
(\Delta_{g_k}-\bar c) E^{(3)} + 
n(n+2) \bar \xi^{ \frac 4{n-2} } 
E^{(3)}=O(M_k^{ -\frac {12}{n-2} })(1+r)^{6-n}, 
\quad 0<r< M_k^{ \frac 2{n-2} }, 
\label{E2-1} 
\end{equation} 
where 
$$ \bar \xi^{\frac{4}{n-2}}(y) 
=\int_0^1(tv_k+(1-t)(U+F^{(3)}))^{\frac 4{n-2}}dt.$$

Arguing as on page 212 of \cite{LiZhang1} , we see that the operator 
$\Delta_{g_k}-\bar c  + n(n+2) \bar \xi^{ \frac 4{n-2} } $ satisfies 
the maximum principle over $R_1<|y|<\delta' M_k^{ \frac 2{n-2} }$ 
for some constants $R_1, \delta'>0$ which are independent of $k$. 
For $C_{10}$ large, but independent of $k$, we see, using 
(\ref{E2-1}), (\ref{nov19e6}) and (\ref{21-new9}), that 
$$ 
\left( \Delta_{g_k}-\bar c 
 + 
n(n+2) \bar \xi^{ \frac 4{n-2} } 
\right)\left( E^{(3)}-f\right)\ge 0, 
\qquad R_1<|y|<\delta' M_k^{ \frac 2{n-2} }, 
$$ 
$$ 
\left( E^{(3)}-f\right)(y)<0\qquad 
\mbox{on}\ \{r=R_1\}\cup \{r= \delta' M_k^{ \frac 2{n-2} }\}, 
$$ 
where $f(r):= C_{10} M_k^{ -\frac {12}{ n-2} }r^{8-n}.$ Thus, in 
view of (\ref{nov19e6}), estimate (\ref{2omega82}) for $l=0$ follows 
from the maximum principle. The estimate for $l=1,2$ can then be 
deduced from the equation satisfied by $E^{(3)}$ using elliptic 
estimates. $\Box$

Recall the Pohozaev type identity (102) in 
\cite{LiZhang1}, with 
$R_k'=\delta' M_k^{ \frac 2{n-2} }$, 
\begin{equation}
I_1[v_k]+I_2[v_k]+I_3[v_k]+I_4[v_k]=I_5[v_k], 
\label{poho}
\end{equation} 
$$ 
I_1[v_k]=\int_{|y|\le R_k'} 
(-\bar b_i \partial_i v_k-\bar d_{ij}\partial_{ij}v_k) 
(\nabla v_k\cdot y+\frac {n-2}2 v_k), 
$$ 
$$ 
I_2[v_k]=-\frac {c(n)}2 
M_k^{ -\frac 4{n-2} } 
\int_{|y|\le R_k'} 
\bigg\{ (M_k^{ -\frac 2{n-2} }y)\cdot \nabla R(M_k^{ -\frac 2{n-2} }y) 
+2R(M_k^{ -\frac 2{n-2} }y)\bigg\} 
v_k^2(y), 
$$ 
$$ 
I_3[v_k]=  \frac {c(n)}2 
M_k^{ -\frac 4{n-2} }   R_k' 
 \int_{  |y|=R_k' } 
R(M_k^{ -\frac 2{n-2} }y)v_k^2(y), 
$$ 
$$ 
I_4[v_k]= -\frac{  (n-2)^2}2 R_k' 
 \int_{  |y|=R_k' } 
v_k(y)^{  \frac {2n}{n-2}}, 
$$ 
$$ 
I_5[v_k]= 
\int_{ |y|=R_k' } 
\bigg\{  (|\frac{\partial v_k}{ \partial \nu}|^2 
-\frac 12|\nabla v_k|^2) R_k' 
+\frac{n-2}2 
v_k \frac{\partial v_k}{ \partial \nu}\bigg\} 
= O(M_k^{ -2}). 
$$

Write 
$$\nabla v_k\cdot y+\frac{n-2}2v_k=\tilde U+\tilde F^{(3)}+ 
\tilde E^{(3)}$$ 
where $$\tilde  U=\nabla U\cdot y+\frac{n-2}2U, 
\quad \tilde F^{(3)}= \nabla  F^{(3)}\cdot y+ \frac{n-2}2 
  F^{(3)}. 
 $$ 
 Clearly, 
\begin{equation} 
|\nabla ^l \tilde F^{(3)}|=O(M_k^{-\frac{8}{n-2}})(1+r)^{6-n-l}, 
\qquad l=0,1,2, 
\label{2omega83} 
\end{equation} 
and 
\begin{equation} 
|\nabla ^l \tilde E^{(3)}|=O(M_k^{-\frac{12}{n-2}}|y|^{8-n-l}), 
\qquad l=0,1,2. 
\label{2omega84} 
\end{equation} 
 
With these and 
(\ref{dec9e1}),(\ref{2omega82}),(\ref{2omega83}),(\ref{2omega84}), 
 we have 
\begin{eqnarray*} 
I_1[v_k]&=&\int_{B(0, M_k^{\frac{2}{n-2}})} 
(-\bar b_i\partial_i-\bar d_{ij}\partial_{ij})(v_k-U) 
(\nabla v_k\cdot y+\frac{n-2}2v_k)dy\\ 
&=&\int_{B(0, M_k^{\frac{2}{n-2}})} 
(\Delta-\Delta_{g_k})\bigg (F^{(3)}+E^{(3)}\bigg ) 
\cdot (\tilde U+\tilde F^{(3)}+\tilde E^{(3)})\\ 
&=& 
\int_{B(0, M_k^{\frac{2}{n-2}})} 
(\Delta-\Delta_{g_k})F^{(3)}\tilde U 
+O(M_k^{-2}). 
\end{eqnarray*} 
Since $\tilde  U$ is radially symmetric, $\int_{B(0, 
M_k^{\frac{2}{n-2}})} (\Delta-\Delta_{g_k})F^{(3)}\tilde U=0.$ Thus 
$ I_1[v_k]=O(M_k^{-2}), $ and with notation in (\ref{15-4new}), 
\begin{eqnarray*} 
I_2[v_k]&=&-\frac{c(n)}2 
\sum_{l=2}^{7}\sum_{|\alpha |=l}\int_{|y|\le R_k'} 
\bigg \{(\frac{l+2}{\alpha!}) 
\partial_{\alpha}Ry^{\alpha}M_k^{-\frac{4+2l}{n-2}} 
\bigg \}v_k^2 
+O(M_k^{-2})\\ 
&=&-\frac{c(n)}2\int_{|y|\le R_k'}\bigg \{ 
\sum_{l=2}^{7}\sum_{p=1}^{I_l}(l+2)\tilde R_{lp}e_pr^l 
M_k^{-\frac{4+2l}{n-2}}+\sum_{s=0}^{2}(2s+4)\bar R^{(2s+2)} 
M_k^{-\frac{8+4s}{n-2}}r^{2s+2}\bigg \}\\ 
&&\quad \cdot \bigg (U^2+2U(F^{(3)}+E^{(3)})+(F^{(3)}+E^{(3)})^2 
\bigg )+O(M_k^{-2}). 
\end{eqnarray*} 
Using Lemma \ref{mar19prop2} and (\ref{dec10e5}), and for 
small $\delta'$, we have 
\begin{eqnarray} 
&& -\frac{c(n)}2\int_{|y|\le R_k'}\bigg \{ 
\sum_{s=0}^{2}(2s+4)\bar R^{(2s+2)} 
M_k^{-\frac{8+4s}{n-2}}r^{2s+2}\bigg \} 
 \bigg (U^2+2U(F^{(3)}+E^{(3)})+(F^{(3)}+E^{(3)})^2 
\bigg )\nonumber 
\\ 
&=& 
- 
\frac{c(n)}2\int_{|y|\le R_k'}\sum_{s=0}^{2} 
(2s+4)\bar R^{(2s+2)}r^{2s+2}U^2 
[ 1+O(\delta')+\circ(1)] 
M_k^{-\frac{8+4s}{n-2}} \nonumber \\ 
&\ge &c_1(n)|W|^2M_k^{-\frac{8}{n-2}}+c_2(n)|\nabla R_{abcd}|^2M_k^{-\frac{12}{n-2}} 
 \nonumber \\ 
&+&4c(n)\bigg (\epsilon |\nabla^2R_{abcd}|^2-\frac{1}{8(n+4)(n+2)n} 
(\frac{n-8}{n-2}-\frac{49}{20n^2}+\epsilon')|\nabla^2 R|^2\bigg ) 
\nonumber \\ 
&&\cdot \bigg ( \int_{|y|\le R_k'}r^6U^2dyM_k^{-\frac{16}{n-2}}\bigg ) 
+\circ(1)\left( |W|^2 M_k^{-\frac{8}{n-2}}+ 
|\nabla R_{abcd}|^2 M_k^{-\frac{12}{n-2}}\right). 
\label{dec17e2} 
\end{eqnarray} 
where $c_1(n),c_2(n), \epsilon$ and $\epsilon'$ are some positive 
constants depending only on $n$. $\epsilon$ and $\epsilon'$ are 
sufficiently small. Also we observe that 
$$-\frac{c(n)}2\int_{|y|\le R_k'}\bigg ( 
\sum_{l=2}^{7}\sum_{p=1}^{I_l}(l+2)\tilde R_{lp}e_pr^l 
M_k^{-\frac{4+2l}{n-2}}\bigg )U^2=0.$$ 
 
We only need to deal with 
 
\begin{eqnarray*} 
&&-\frac{c(n)}2\int_{|y|\le R_k'}\bigg ( 
\sum_{l=2}^{7}\sum_{p=1}^{I_l}(l+2)\tilde R_{lp}e_pr^l 
M_k^{-\frac{4+2l}{n-2}}+\sum_{s=0}^{2}(2s+4)\bar R^{(2s+2)} 
M_k^{-\frac{8+4s}{n-2}}r^{2s+2}\bigg )\\ 
&&\cdot \bigg (2U(F^{(3)}+E^{(3)})+(F^{(3)}+E^{(3)})^2 
\bigg ). 
\end{eqnarray*} 
 
By previous estimates 
 
$$\int_{|y|\le R_k'}\bigg (\sum_{l=3}^{7}\sum_{p=1}^{I_l}(l+2)\tilde R_{lp}e_pr^l 
M_k^{-\frac{4+2l}{n-2}}\bigg ) 
\bigg (2U(F^{(3)}+E^{(3)})+(F^{(3)}+E^{(3)})^2 \bigg )=O(M_k^{-2}).$$ 
 
$$\int_{|y|\le R_k'}\bigg (\sum_{p=1}^{I_2}\tilde R_{2p}e_pr^2 
M_k^{-\frac{8}{n-2}}\bigg )\bigg (2UE^{(3)}+(F^{(3)}+E^{(3)})^2\bigg 
) =O(M_k^{-2}).$$ 
 
Finally $-\frac{c(n)}2\int_{|y|\le R_k'}\sum_{l=2}^{7}\sum_{p=1}^{I_l} 
(l+2)\tilde R_{lp}e_pr^lM_k^{-\frac{4+2l}{n-2}}2UF^{(3)}$ contributes another important 
term. 
 
\begin{eqnarray*} 
&&-\frac{c(n)}2\int_{|y|\le R_k'}\sum_{l=2}^{7}\sum_{p=1}^{I_l} 
(l+2)\tilde R_{lp}e_pr^lM_k^{-\frac{4+2l}{n-2}}2UF^{(3)}\\ 
&=&-c(n)\int_{|y|\le R_k'} 
\bigg (\sum_{p=1}^{I_2}4\tilde R_{2p}e_p r^2M_k^{-\frac{8}{n-2}} 
+\sum_{l=3}^{7}\sum_{p=1}^{I_l} 
(l+2)\tilde R_{lp}e_pr^lM_k^{-\frac{4+2l}{n-2}}\bigg )\\ 
&&U(F^{(2)}+F^{(3)}-F^{(2)}). 
\end{eqnarray*} 
 
Using the fact the eigenfunction corresponding to $l=2$ are orthogonal to those 
corresponding to $l=3$, we have 
\begin{eqnarray} 
&&-\frac{c(n)}2\int_{|y|\le R_k'}\sum_{l=2}^{7}\sum_{p=1}^{I_l} 
(l+2)\tilde R_{lp}e_pr^lM_k^{-\frac{4+2l}{n-2}}2UF^{(3)} \nonumber \\ 
&=&-4c(n)\int_{|y|\le R_k'} 
\widetilde R^{(2)}(\theta) 
r^2UF^{(2)}M_k^{-\frac 8{n-2}}+O(M_k^{-2}) 
+\circ(1)M_k^{-\frac 8{n-2}}|W|^2
\nonumber \\ 
&=&\frac{2c(n)^2}{n(n+2)}\bigg (\sum_{i<j} 
2(\partial_{ij}R)^2+\sum_i(\partial_{ii}R)^2 \bigg) 
\int_{|y|\le R_k'}r^2Uf_{2}dyM_k^{-\frac{16}{n-2}}+O(M_k^{-2}) \nonumber \\ 
&&+\circ(1)M_k^{-\frac 8{n-2}}|W|^2
\nonumber\\
&\ge &\frac{2c(n)^2}{n(n+2)}|\nabla^2R|^2\int_{|y|\le 
R_k'}r^{2}Uf_{2}dy M_k^{-\frac{16}{n-2}}+O(M_k^{-2}) 
+\circ(1)M_k^{-\frac 8{n-2}}|W|^2
 \label{dec17e3} 
\end{eqnarray} 
The Pohozaev type identity (\ref{poho})
yields, in view of 
 (\ref{dec17e2}), (\ref{dec17e3}) and the lower bound of $f_2$ in 
(\ref{nov19e1}),
that
\begin{eqnarray*} 
|W|^2M_k^{-\frac{8}{n-2}}+|\nabla R_{abcd}|^2M_k^{-\frac{12}{n-2}} 
+|\nabla^2 R_{abcd}|^2M_k^{-\frac{16}{n-2}}\log M_k&=&O(M_k^{-2}),\quad n=10.\\ 
\\ 
|W|^2M_k^{-\frac{8}{n-2}}+|\nabla R_{abcd}|^2M_k^{-\frac{12}{n-2}} 
+|\nabla^2 R_{abcd}|^2M_k^{-\frac{16}{n-2}}&=&O(M_k^{-2}),\quad n=11. 
\end{eqnarray*} 
Thus we have proved 
(\ref{W1}) and (\ref{W2}).
Estimates (\ref{W7}) and (\ref{W8}) follow from
(\ref{2omega82}). 
Theorem \ref{thm1} 
is established. 
 $\Box$

\section{Appendix A: Some curvature inequalities
in
 conformal normal coordinates}

\subsection{The inequality for  $\bar R^{(6)}$} 
 
In this subsection we prove the following 
two lemmas.
\begin{lem} 
If $|W(0)|=|\nabla W(0)|=0$, then we have, in  conformal normal coordinates 
centered at $0$, 
\begin{eqnarray} 
\bar R^{(6)}&=&-\frac{R_{p_1p_2p_3p_4,p_5p_6}R_{p_1p_2p_3p_4,p_5p_6}} 
{40(n+4)(n+2)n}-\frac{R_{p_1p_2,p_3p_4}(R_{p_1p_2,p_3p_4}+R_{p_3p_4,p_1p_2})} 
{8(n+4)(n+2)n} \nonumber \\ 
&& +\frac{\sum_{p_1p_2}(\partial_{p_1p_2}R)^2}{8(n+4)(n+2)n},\quad 
\mbox{at}\quad 0, 
\label{dec7e4} 
\end{eqnarray} 
where repeated indices mean summation, and $R_{ijkl, pq}$ denotes
convariant derivatives of $R_{ijkl}$.
\label{mar19prop1} 
\end{lem} 
 
\begin{lem} 
For some small $\epsilon=\epsilon(n)>0$, we have,
in  conformal normal coordinates
centered at $0$,
\begin{eqnarray} 
&&\bar R^{(6)}<-\epsilon R_{p_1p_2p_3p_4,p_5p_6}R_{p_1p_2p_3p_4,p_5p_6} 
\nonumber \\ 
&&+\frac{1}{8(n+4)(n+2)n}(\frac{n-8}{n-2}- 
\frac{49}{20n^2}+\epsilon)\sum_{p_1p_2}(\partial_{p_1p_2}R)^2 
+O(|\nabla R_{abcd}|^2+|W|^2). \label{dec10e8} 
\end{eqnarray} 
\label{mar19prop2} 
\end{lem} 
 
We first assume Lemma \ref{mar19prop1} and give the proof of 
Lemma \ref{mar19prop2}. 
 
\noindent{\bf Proof of Lemma \ref{mar19prop2}:} 

It was proved by Hebey and Vaugon in  \cite{HV} that,
if  $|W(0)|=|\nabla
W(0)|=0$, then, in  conformal normal coordinates
centered at $0$,
\begin{equation} 
R_{p_1p_2,p_3p_4}(R_{p_1p_2,p_3p_4}+R_{p_3p_4,p_1p_2})\ge 
\frac{6}{n-2}\sum_{p_1p_2}(\partial_{p_1p_2}R)^2,\quad \mbox{at}\quad 0. 
\label{dec7e2} 
\end{equation} 
We also need the following inequality under the same assumption: 
\begin{equation} 
R_{p_1p_2p_3p_4,p_5p_6}R_{p_1p_2p_3p_4,p_5p_6}\ge \frac{49}{4n^2} 
\sum_{p_1p_2}(\partial_{p_1p_2}R)^2,\quad \mbox{at}\quad 0. 
\label{dec7e3} 
\end{equation} 
Note that (\ref{dec7e3}) with $\frac{49}{4n^2}$ replaced by 
$\frac{9}{n^2}$ was established  in \cite{HV}. 
This weaker version leads to an
  inequality weaker  than (\ref{mar19prop2}), which is nevertheless 
enough for applications in this paper.
 To  prove 
(\ref{dec7e3}), we consider 
$$\|R_{ikmj,pq}-\alpha R_{,ij}\delta_{kp}\delta_{mq}\|^2>0.$$ 
Namely,
$$|\nabla_{pq}R_{ikmj}|^2-2\alpha R_{ikmj,km}R_{,ij} 
+\alpha^2n^2R_{,ij}R_{,ij}\ge 0.$$  By the
 second Bianchi identity,
$R_{ikmj,k}=R_{im,j}-R_{ij,m}$. So 
$$R_{ikmj,km}=R_{im,jm}-R_{ij,mm}=\frac 12R_{,ij}+3R_{,ij}=\frac 72R_{,ij}$$ 

It follows that 
$$|\nabla_{pq}R_{ikmj}|^2+(\alpha^2n^2-7\alpha)R_{,ij}R_{,ij}\ge 0.$$ 
Inequality (\ref{dec7e3}) follows from the above by taking
$\alpha=\frac 7{2n^2}$.

By (\ref{dec7e4}), (\ref{dec7e2}) and (\ref{dec7e3}) we have 
$$\bar R^{(6)}<\frac{1}{8(n+4)(n+2)n}(\frac{n-8}{n-2}- 
\frac{49}{20n^2})\sum_{p_1p_2}(\partial_{p_1p_2}R)^2,\quad 
\mbox{at}\quad 0.$$ Then (\ref{dec7e8}) holds for some small 
$\epsilon(n)>0$ under the assumption $|W(0)|=|\nabla W(0)|=0$. In 
general if we do not assume $|W(0)|=|\nabla W(0)|=0$,  all the extra 
terms can be estimated by Cauchy's inequality, and  we  obtain 
(\ref{dec10e8}).  Lemma \ref{mar19prop2} is established. $\Box$ 
 
\medskip 
 
\noindent{\bf Proof of Lemma \ref{mar19prop1}:} 

It was proved in \cite{HV} that if  $|W(0)|=|\nabla W(0)|=0$, 
then, in conformal normal coordinates centered at  $0$,
\begin{eqnarray*} 
&&C(2,2)\mbox{Sym}_{p_1..p_6}R_{,p_1..p_6}+864R_{p_1p_2p_3p_4,p_5p_6} 
R_{p_1p_2p_3p_4,p_5p_6}\\ 
&&+4320R_{p_1p_2,p_3p_4}(R_{p_1p_2,p_3p_4}+R_{p_3p_4,p_1p_2}) 
-4320\sum_{p_1p_2}(\partial_{p_1p_2}R)^2=0, 
\end{eqnarray*} 
where $C(2,2)$ is the complete contraction: 
$$C(2,2)=\sum_{p_1=p_2=1}^n\sum_{p_3=p_4=1}^n\sum_{p_5=p_6=1}^n.$$ 
 Since we work in conformal normal coordinates and since
$|W(0)|=|\nabla W(0)|=0$,
$$\mbox{Sym}_{p_1..p_6}R_{,p_1..p_6}=\mbox{Sym}_{p_1..p_6}\partial_{p_1..p_6}R
-\frac{36}5\mbox{Sym}_{p_1..p_6}R_{,\nu p_1}R_{\nu
p_2p_3p_4,p_5p_6}
=\mbox{Sym}_{p_1..p_6}\partial_{p_1..p_6}R,
$$
where, for the second equality, we have used 
 the skew-symmetry of
$R_{abcd}$). 

Thus we have
$$C(2,2)\mbox{Sym}_{p_1..p_6}R_{,p_1..p_6}=720\Delta^3R(0),$$
where $\Delta$ denotes the flat Laplacian.
 
Therefore 
\begin{eqnarray} 
&&\Delta^3R(0)+\frac 65R_{p_1p_2p_3p_4,p_5p_6}R_{p_1p_2p_3p_4,p_5p_6} 
+6R_{p_1p_2,p_3p_4}(R_{p_1p_2,p_3p_4}+R_{p_3p_4,p_1p_2}) \nonumber \\ 
&&-6\sum_{p_1p_2}(\partial_{p_1p_2}R)^2=0. 
\label{dec13e1} 
\end{eqnarray} 
By some standard computations,
\begin{equation} 
\bar R^{(6)}=\frac 1{|S^{n-1}|}\int_{S^{n-1}}\sum_{|\alpha |=6} 
\frac{\partial_{\alpha}R}{\alpha !}\theta^{\alpha}
=\frac{\Delta^3R(0)}{48(n+4)(n+2)n}. 
\label{dec13e4} 
\end{equation} 
It follows from (\ref{dec13e1}) and (\ref{dec13e4}) that 
\begin{eqnarray*} 
\bar R^{(6)}&=&-\frac{R_{p_1p_2p_3p_4,p_5p_6}R_{p_1p_2p_3p_4,p_5p_6}} 
{40(n+4)(n+2)n}-\frac{R_{p_1p_2,p_3p_4}(R_{p_1p_2,p_3p_4}+R_{p_3p_4,p_1p_2})} 
{8(n+4)(n+2)n}\\ 
&& +\frac{\sum_{p_1p_2}(\partial_{p_1p_2}R)^2}{8(n+4)(n+2)n} 
\end{eqnarray*} 
at $0$ where $|W(0)|=|\nabla W(0)|=0$ is assumed. 
Lemma \ref{mar19prop1} is established. $\Box$

\subsection{Proof
of Lemma \ref{lemjune27}} 

The following fact is elementary: 
Let $k\ge 1$ be an integer, $j\in [1,n]$ be a fixed integer, then 
\begin{equation}
\int_{S^{n-1}}\partial_{p_1..p_{2k+1}}R(0)x^{p_1}
\cdots x^{p_{2k+1}}\cdot x^jdx= 
C(n,k)\partial_j(\Delta^kR)(0)
\label{fact1}
\end{equation}
 where 
$$
C(n,k)=\frac{(2k+1)!|S^{n-1}|}{(2k+n)2^kk!\prod_{i=0}^{k-1}(n+2i)}.$$

\noindent{\bf Proof of Lemma \ref{lemjune27}:} 
In
conformal normal coordinates,
$$\mbox{Sym}_{p_1\cdots p_{2k+3}}R_{p_1p_2,p_3, \cdots, p_{2k+3}}=0,\quad 
\omega+2\le 2k+3\le 2\omega+3,$$ if $|\nabla^iR_{abcd}(0)|=0$ for 
$0\le i\le \omega-1$. See \cite{HV}. After contraction this implies 
$$\partial_j(\Delta^k)R(0)=0 \quad j=1,
\cdots, n,\quad \mbox{if}\quad 
|\nabla^iR_{abcd}(0)|=0\quad \mbox{for}\quad 0\le i\le \omega-1.$$ 
For $ n=10,  11$, we only need to discuss $k=1$, i.e, we have 
$\partial_j(\Delta R)(0)=0$ if $|W(0)|=0$. In general we have 
$$\partial_j(\Delta R)(0)=O(|W|).$$ 
This and (\ref{fact1}) imply
$$
\int_{\Bbb S^{n-1} } \tilde R^{(3)}(\theta)\theta^j=
O(|W|),\qquad 1\le j\le n.
$$
On the other hand, it is clear that
$$
\int_{\Bbb S^{n-1} } \tilde R^{(3)}(\theta)=
0, \quad 
\int_{\Bbb S^{n-1} } \tilde R^{(3)}(\theta)\theta^i\theta^j=0\quad
1\le i,j\le n.
$$
Lemma \ref{lemjune27} follows from the above.
$\Box$

\section{Appendix B: Some estimates on an ODE}

\begin{prop} 
\label{compprop4} 
Let  $n\ge 3$ be an integer, $\delta_0\ge n$ be a constant, let $\hat H(r)\in C^0(0,\infty)$ 
satisfy, for some positive constants $C,\beta$ and $\alpha>2$, 
$\delta_0+(\alpha -2)(n-\alpha )>0$, 
$$ 
0\le \hat H(r)\le Cr^{\beta}(1+r)^{-\beta-\alpha}, \quad 0<r<\infty. 
$$ 
Then for any constant $p$ satisfying 
$$ 
0<p\le \beta+2,\quad p(p+n-2)<\delta_0, 
$$ 
there exists a unique $a(r)\in C^2(0,\infty)$ verifying 
\begin{equation} 
\left\{\begin{array}{ll} 
Ta(r):=a''(r)+\frac{n-1}ra'(r)+(n(n+2)U(r)^{\frac{4}{n-2}}-\frac{\delta_0}{r^2})a(r)= 
-\hat H(r),\quad 0<r<\infty, \\ 
\\ 
\lim_{r\to 0}a(r)=\lim_{r\to \infty}a(r)=0. 
\end{array} 
\right. 
\label{comp1} 
\end{equation} 
Moreover, for some positive constant $C_0$ depending only on $n$, 
$\delta_0$, $\alpha$, $\beta$, $p$ and $C$, 
$$ 
0\le a(r)\le C_0r^p(1+r)^{-p+2-\alpha},\quad 0<r<\infty. 
$$ 
\end{prop} 
 
\begin{lem} Let  $n\ge 3$ be an integer, 
 $\delta_0\ge n$ be a constant, and let $\hat H(r)$ be a non-negative 
function in $C^0(0,\infty)$.  Then 
for any $0<\epsilon<R$, there exists a unique solution 
$a_{\epsilon,R}\in C^2[\epsilon, R]$ to 
\begin{equation} 
\left\{\begin{array}{ll} 
a_{\epsilon,R}''(r)+\frac{n-1}ra_{\epsilon,R}'(r) 
+(n(n+2)U(r)^{\frac{4}{n-2}}-\frac{\delta_0}{r^2})a_{\epsilon,R}(r)= 
-\hat H(r),\ \ \epsilon<r<R,\\ 
\\ 
a_{\epsilon,R}(\epsilon)=a_{\epsilon,R}(R)=0, 
\end{array} 
\right. 
\label{dec5e1} 
\end{equation} 
Moreover, $a_{\epsilon, R}\ge 0$ on $[\epsilon, R]$. 
\label{lemB.1} 
\end{lem} 
 
\noindent{\bf Proof.}\ 
Let $(S^n,g_0)$ be the standard sphere. It is known that in the 
stereographic projection coordinates 
$$g_0=\sum_{i=1}^{n+1}dx_i^2=u_1(y)^{\frac{4}{n-2}}dy^2$$ 
where 
$$u_1(y)=(\frac{2}{1+|y|^2})^{\frac{n-2}{2}}=2^{\frac{n-2}{2}}U.$$ 
Also we know that 
$$L_{g_0}(\phi)=(\Delta_{g_0}-\frac{n(n-2)}{4})\phi=u_1^{-\frac{n+2}{n-2}} 
\Delta (u_1\phi).$$ 
If we let $\phi=a_{\epsilon,R}/u_1$, we can rewrite (\ref{dec5e1}) as 
 
$$(\Delta_{g_0}\phi-\frac{n(n-2)}{4}\phi)u_1^{\frac{n+2}{n-2}} 
=-(n(n+2)U^{\frac{4}{n-2}}-\frac{\delta_0}{r^2})a_{\epsilon,R}-\hat H(r).$$ 
After simplification, we have 
$$\Delta_{g_0}\phi+(n-\frac{\delta_0(1+r^2)^2}{4r^2})\phi=-u_1^{-\frac{n+2}{n-2}} 
\hat H(r).$$ 
Since $\delta_0\ge n$ 
\begin{equation} 
n<\frac{\delta_0(1+r^2)^2}{4r^2},\quad \mbox{for}\quad \epsilon\le r\le R, 
\quad r\neq 1. 
\label{6-1} 
\end{equation} 
So the existence and the uniqueness of $\phi$ as well as 
$a_{\epsilon,R}$ are proved. $\Box$ 
 
\medskip 
 
\noindent{\bf Proof of Proposition \ref{compprop4}.}\ Clearly for 
some $R_1>1$, 
$$\frac{\delta_0}{r^2}\ge \frac{8n}{r^4},\quad \mbox{ for } r\ge R_1/2,\qquad \mbox{and} 
\quad \frac{\delta_0}{r^2}\ge 8n,\quad \mbox{ for } 0<r\le 
\frac{2}{R_1}.$$ Fix a $W\in C^0(0,\infty)$ satisfying 
\begin{eqnarray*} 
W(r)=\frac{\delta_0}{r^2},\quad && \frac 2{R_1}<r<\frac{R_1}2,\\ 
\frac{\delta_0}{r^2}\ge W(r)\ge 4n,\quad && \frac 1{R_1}<r<\frac 2{R_1},\\ 
W(r)=4n,\quad &&0<r\le \frac 1{R_1},\\ 
\frac{\delta_0}{r^2}\ge W(r)\ge \frac{4n}{r^4},\quad && 
\frac{R_1}2<r\le R_1,\\ 
W(r)=\frac{4n}{r^4},\quad &&r>R_1. 
\end{eqnarray*} 
 
Clearly $W(r)\le \delta_0/r^2$ and 
$$\frac{(1+r^2)^2}{4}W(r)>n,\quad \forall r>0,\quad r\neq 1.$$ 
With this fact, the first eigenvalue of 
$-\Delta_{g_0}+\bigg (\frac{(1+r^2)^2}{4}W(r)-n\bigg )$ is positive on 
$S^n$ and the potential 
$\bigg (\frac{(1+r^2)^2}{4}W(r)-n\bigg )$ is in $C^0(S^n)$. Since 
$\alpha >2$, $u_1^{-\frac{n+2}{n-2}}\hat H(r)\in L^q(S^n)$ for some 
$q>1$. Let $\phi_1\in W^{2,q}(S^n)$ be the solution of 
$$\Delta_{g_0}\phi_1+\bigg (n-\frac{(1+r^2)^2}{4}W(r)\bigg )\phi_1 
=-u_1^{-\frac{n+2}{n-2}}\hat H(r).$$ 
By the symmetry of the data, the uniqueness of the solution, $\phi_1$ 
depends only on $r$. Since both $W$ and $\hat H$ are continuous for 
$0<r<\infty$, $\phi_1$ is $C^2$ in $0<r<\infty$. By the maximum principle, 
$\phi_1\ge 0$. Let 
$$a_1(r)=\phi_1(r)u_1(r),\quad 0<r<\infty.$$ 
Then $a_1\in C^2(0,\infty)$, $a_1(r)\ge 0,$ and 
\begin{eqnarray*} 
Ta_1(r):&=&a_1''(r)+\frac{n-1}ra_1'(r)+ 
[n(n+2)U(r)^{\frac 4{n-2}}-\frac{\delta_0}{r^2}]a_1(r)\\ 
&\le & a_1''(r)+\frac{n-1}ra_1'(r)+ 
[n(n+2)U(r)^{\frac 4{n-2}}-W(r)]a_1(r)\\ 
&=&(\Delta_{g_0}\phi_1-\frac{n(n-2)}{4}\phi_1)u_1^{\frac{n+2}{n-2}} 
-[n(n+2)U(r)^{\frac 4{n-2}}-W(r)]\phi_1u_1\\ 
&=&\bigg (\Delta_{g_0}\phi_1+(n-\frac{(1+r^2)^2}{4}W(r))\phi_1\bigg ) 
u_1^{\frac{n+2}{n-2}}=-\hat H_1(r). 
\end{eqnarray*} 
 
So $a_1(r)$ is a supersolution. A calculation gives 
\begin{eqnarray*} 
T(r^p)&=&-\bigg (\delta_0-p(p+n-2)+O(r)\bigg )r^{p-2},\quad \mbox{as}\quad r\to 0,\\ 
T(r^{2-\alpha})&=&-\bigg ( 
\delta_0+(\alpha-2)(n-\alpha)+O(1/r)\bigg )r^{-\alpha},\quad \mbox{as} 
\quad r\to \infty. 
\end{eqnarray*} 
Since both $\delta_0-p(p+n-2)$ and $\delta_0+(\alpha-2)(n-\alpha)$ are positive, 
there exists $R_2>1$ such that 
\begin{eqnarray*} 
T(\gamma r^p)&\le &-\hat H(r),\quad \mbox{for}\quad 0<r\le \frac 1{R_2},\\ 
T(\gamma r^{2-\alpha})&\le &-\hat H(r),\quad \mbox{for}\quad r\ge R_2 
\end{eqnarray*} 
for some $\gamma>1$. Choose $\gamma$ larger if necessary such that 
$$\gamma(\frac 1{R_2})^p>a_1(\frac 1{R_2}),\quad \gamma (R_2)^{2-\alpha}>a_1(R_2).$$ 
Define 
$$\bar a(r)=\left\{\begin{array}{ll} 
\min\{\gamma r^p,a_1(r)\},&\quad 0<r<\frac 1{R_2},\\ \\ 
a_1(r),&\quad \frac 1{R_2}\le r\le R_2,\\ \\ 
\min\{\gamma r^{2-\alpha},a_1(r)\},&\quad r>R_2. 
\end{array} 
\right. 
$$ 
Then $\bar a(r)$ is a continuous supersolution to $T\bar a(r)=-\hat H(r)$ in 
$(0,\infty)$. Therefore for any $0<\epsilon<R<\infty$, then solution of (\ref{dec5e1}) 
satisfies 
$$0\le a_{\epsilon,R}(r)\le \bar a(r),\quad \forall 0<r<\infty.$$ 
Let $\epsilon\to 0$ and $R\to \infty$ along a subsequence, $a_{\epsilon,R}(r)$ tends 
to $a(r)$ in $C^{1,\lambda}_{loc}(0,\infty)$ for $0<\lambda<1$, which satisfies 
(\ref{comp1}) in the weak sense. Since $\hat H\in C^0(0,\infty)$, we know that 
$a\in C^2(0,\infty)$. 
 
Now we prove the uniqueness of the solution of (\ref{comp1}). Let $a(r)$ and $b(r)$ be 
two solutions of (\ref{comp1}), then their difference verifies the homogeneous 
equation 
$$T(a-b)\equiv 0,\quad \mbox{in}\quad (0,\infty).$$ 
By Theorem 8.1 in \cite{CLev}, the homogeneous equation has two linearly 
independent solutions $a_{+}(r)$ and $a_{-}(r)$ with the asymptotic behavior 
$$\lim_{r\to \infty}\frac{a_+(r)}{r^{\lambda_1}}=\lim_{r\to \infty} 
\frac{a_-(r)}{r^{\lambda_2}}=1,$$ 
where $\lambda_1$ and $\lambda_2$ are the two solutions of 
$\lambda^2+(n-2)\lambda-\delta_0=0$ such that $\lambda_1>0$ and 
$\lambda_2<2-n$. Since $(a-b)(r)=C_1a_++C_2a_-$ for some constants $C_1$ and $C_2$, and 
since $(a-b)(r)\to 0$ as $r\to \infty$, we must have $C_1=0$ and therefore 
\begin{equation} 
\lim_{r\to \infty}r^{n-2}(a-b)(r)=0. 
\label{D10-1} 
\end{equation} 
Since $T(a-b)=0$ corresponds to 
$$(\Delta_{g_0}+[n-\frac{\delta_0(1+r^2)^2}{4r^2}]) 
\frac{(a-b)(r)}{u_1}=0,\quad \mbox{on}\quad S^n\setminus \{P,N\},$$ 
where $N,P$ are the south pole and the north pole, respectively. We 
know from (\ref{D10-1}) that $\frac{(a-b)}{u_1}(p)\to 0$ as $p\to 
\{N,P\}$, and in view of (\ref{6-1}), we can apply the maximum 
principle to conclude $a-b\equiv 0$. Proposition \ref{compprop4} is 
established. $\Box$

\section{Appendix C: Two useful lower bounds }

\begin{prop} 
\label{propmar11} 
For $n\ge 10$, there exists a unique 
$f_2\in C^\infty( (0,\infty))$ satisfying 
\begin{equation} 
\left\{\begin{array}{ll} 
f_2''(r)+\frac{n-1}rf_2'(r)+(n(n+2)U^{\frac{4}{n-2}}-\frac{2n}{r^2})f_2 
=-r^2U,\quad 0<r<\infty,\\ 
\lim_{r\to 0}f_2(r)=\lim_{r\to \infty}f_2(r)=0. 
\end{array} 
\right. 
\label{nov19e4} 
\end{equation} 
Moreover, for some 
universal positive constant $C$, 
\begin{equation} 
\label{nov19e1} 
\frac{U}{6(n-4)}(r^4+\frac{3n-4}{n-2}r^2) 
\le f_2(r)\le Cr^{\frac 32}(1+r)^{\frac 92-n}, 
\qquad 0<r<\infty, 
\end{equation} 
\end{prop}

\begin{prop} 
Let $V_{\lambda}=n(n+2)\int_0^1(tU+(1-t)U^{\lambda})^{\frac 
4{n-2}}dt$.  Then there exists a unique
  $f_{2,\lambda}\in C^\infty (0,\infty)$
satisfying 
$$ 
\left\{\begin{array}{ll} 
f_{2,\lambda}''(r)+\frac{n-1}rf_{2,\lambda}'(r)+ 
(V_{\lambda}-\frac{2n}{r^2})f_{2,\lambda}(r)=-r^2U^{\lambda}(r) 
(1-(\frac{\lambda}{r})^{8}),\quad r\in (\lambda,\infty)\\ 
\\ 
f_{2,\lambda}(\lambda)=0,\quad \lim_{r\to \infty}f_{2,\lambda}(r)=0. 
\end{array} 
\right. 
$$ 
Moreover for any $\epsilon>0$, there exist $\delta(\epsilon)$ satisfying 
$\delta(\epsilon)\to 0$ as $\epsilon\to 0$, and a universal constant 
$C$ such that for $|\lambda-1|\le \delta(\epsilon)$, 
\begin{equation} 
\frac{1-\epsilon}{6(n-4)}U^{\lambda} 
\bigg (r^4(1-(\frac{\lambda}{r})^8)+\frac{3n-4}{n-2}r^2(1-(\frac{\lambda}{r})^4)\bigg ) 
\le f_{2,\lambda}(r)\le Cr^{6-n} 
\label{mar10e2} 
\end{equation} 
for $\lambda<r<\infty$. 
\end{prop} 
 
\begin{prop} 
For $n\ge 8$, there exists a unique $f_3\in C^\infty (0,\infty)$ 
satisfying 
\begin{equation} 
\left\{\begin{array}{ll} 
f_3''(r)+\frac{n-1}rf_3'(r)+(n(n+2)U^{\frac{4}{n-2}}-\frac{3(n+1)}{r^2})f_3 
=-r^3U,\quad 0<r<\infty,\\ 
\\ 
\lim_{r\to 0}f_3(r)=\lim_{r\to \infty}f_3(r)=0. 
\end{array} 
\right. 
\label{aabb}
\end{equation}
and, for a universal constant $C$ 
$$ 
0\le f_3(r)\le Cr^{\frac 52}(1+r)^{\frac 92-n},\quad r>0. 
$$ 
\end{prop}

The existence, uniqueness, and the 
upper bounds of $f_2$, $f_{2,\lambda}$ and $f_3$ follow
from  Proposition \ref{compprop4}. 
So we only prove the lower bound of $f_2$ and $f_{2,\lambda}$ in this section. 
 
Let $\phi_1(r)=r^4U/(6(n-4))$. Then by elementary computation 
$$ 
\Delta \phi_1+\bigg (n(n+2)U^{\frac 4{n-2}}-\frac{2n}{r^2} \bigg 
)\phi_1>-r^2U. 
$$ 
Consequently 
$$\left\{\begin{array}{ll} 
T(f_2- \phi_1):=\Delta (f_2- \phi_1)+(n(n+2)U^{\frac{4}{n-2}} 
-\frac{2n}{r^2})(f_2- \phi_1)=-g\le 0, 
\quad 0<r<\infty,\\ \\ 
\lim_{r\searrow 0}(f_2-\phi_1)(r)=\lim_{r\to 
\infty}(f_2-\phi_1)(r)=0 
\end{array} 
\right. 
$$ 
where 
$$ 
g(r)=r^2U\bigg (\frac{4(n-2)}{3(n-4)}\frac 1{1+r^2}+\frac{2n}{3(n-4)}r^2U^{\frac 4{n-2}}\bigg ). 
$$ 
 
By Proposition \ref{compprop4}, there exists a positive solution $a_1(r)$ of 
$$\left\{\begin{array}{ll} 
Ta_1(r)=-g(r),\quad 0<r<\infty,\\ 
\\ 
\lim_{r\to 0}a_1(r)=\lim_{r\to \infty}a_1(r)=0. 
\end{array} 
\right. 
$$ 
Since $$\left\{\begin{array}{ll} 
T(f_2-\phi_1-a_1)=0,\quad 0<r<\infty,\\ 
\\ 
\lim_{r\to 0}(f_2-\phi_1-a_1)(r)=\lim_{r\to \infty}(f_2- 
\phi_1-a_1)(r)=0, 
\end{array} 
\right. 
$$ 
we know from the proof of Proposition \ref{compprop4} that $f_2- 
\phi_1-a_1\equiv 0$. So we only need to obtain a lower bound for 
$a_1$. Let 
$$\phi_2(r)=\frac{3n-4}{6(n-4)(n-2)}r^2U,$$ then direct computation 
gives 
$$ 
\Delta \phi_2+(n(n+2)U^{\frac{4}{n-2}}-\frac{2n}{r^2})\phi_2 
=-\frac{2(3n-4)}{3(n-4)}U(\frac{r^2}{1+r^2}-\frac{n}{n-2}r^2U^{\frac{4}{n-2}}). 
$$ 
Then one verifies immediately 
$$g(r) 
> \frac{2(3n-4)}{3(n-4)}U(\frac{r^2}{1+r^2}-\frac{n}{n-2}r^2 
U^{\frac{4}{n-2}}).$$ 
This means 
$$f_2(r)-\frac{1}{6(n-4)}r^4U>\frac{3n-4}{6(n-4)(n-2)}r^2U.$$ 
(\ref{nov19e1}) is established. 
 
\medskip 
 
To prove (\ref{mar10e2}), we still use maximum principle as before, but 
instead of comparing $f_{2,\lambda}$ directly with the right-hand-side in 
(\ref{mar10e2}), we compare $f_{2,\lambda}$ with \\ 
$(1-2\epsilon)(\phi_3+\phi_4)$ where 
\begin{eqnarray*} 
&&\phi_3:=\phi_1-\phi_1^{\lambda}=\frac{1}{6(n-4)}\bigg 
(r^4U^{\lambda}(1-(\frac{\lambda}{r}))^8 
+r^4(U-U^{\lambda})\bigg ),\\ 
&&\phi_4:=\phi_2-\phi_2^{\lambda}=\frac{3n-4}{6(n-4)(n-2)}\bigg 
(r^2U^{\lambda}(1-(\frac{\lambda}{r}))^4 +r^2(U-U^{\lambda})\bigg ). 
\end{eqnarray*} 
 
Our purpose is to show 
\begin{equation} 
f_{2,\lambda}\ge (1-2\epsilon)(\phi_3+\phi_4) \label{mar11e1} 
\end{equation} 
where $\epsilon$ is any fixed small positive constant and $\lambda$ 
is close to $1$ depending on $\epsilon$. Once we have 
(\ref{mar11e1}), (\ref{mar10e2}) follows from (\ref{mar11e1}) and 
the following well known fact: 
 
$$ 
U(r)-U^{\lambda}(r)=(1-\lambda)(1-\frac{\lambda}{r})O(r^{2-n}). 
$$ 
 
Let $T:=\Delta+V_{\lambda}-\frac{2n}{r^2},$ then by elementary 
computation, 
\begin{eqnarray*} 
T\phi_3&=&-r^2U^{\lambda}(1-(\frac{\lambda}{r})^8)+\frac{4(n-2)}{3(n-4)}U^{\lambda} 
(\frac{r^2}{1+r^2}-\frac{\lambda^8}{r^6(1+\lambda^4/r^2)}) \nonumber \\ 
&&+\frac{2n}{3(n-4)}r^4(U^{\lambda})^{\frac{n+2}{n-2}}(1-(\frac{\lambda}{r})^8) 
+\delta(\lambda)(1-\frac{\lambda}r)r^{4-n}. 
\end{eqnarray*} 
 
\begin{eqnarray*} 
T\phi_4&=&-\frac{2(3n-4)}{3(n-4)}U^{\lambda}(1-(\frac{\lambda}{r})^4)+\frac{2(3n-4)}{3(n-4)}U^{\lambda} 
(\frac{1}{1+r^2}-(\frac{\lambda}r)^4\frac 1{1+\lambda^4/r^2}) \nonumber \\ 
&&+\frac{2n(3n-4)}{3(n-4)(n-2)}r^2(U^{\lambda})^{\frac{n+2}{n-2}}(1-(\frac{\lambda}{r})^4) 
+\delta(\lambda)(1-\frac{\lambda}r)r^{2-n}. 
\end{eqnarray*} 
where we use $\delta(\lambda)$ to indicate a function of $\lambda$ 
which tends to 0 as $\lambda\to 1$. 
 
Then one verifies that 
$$T\bigg (f_{2,\lambda}-(1-2\epsilon)(\phi_3 
+\phi_4)\bigg ) < 0,\quad \lambda<r<\infty$$ if $\lambda$ is close 
to $1$ enough. Since we  also have 
$f_{2,\lambda}(\lambda)=\phi_3(\lambda)=\phi_4(\lambda)=0$ and 
$\lim_{r\to \infty}(f_{2,\lambda}-(1-2\epsilon) (\phi_3+\phi_4))=0$, 
we have proved (\ref{mar11e1}) by maximum principle. Proposition 
\ref{propmar11} is established. $\Box$

\end{document}